\documentclass[twocolumn]{autart}    

\usepackage{graphicx}          

 \usepackage{amsfonts,amsmath,amssymb,bm}
\usepackage{color,cases,float}
\usepackage{graphicx}
\usepackage{hyperref}
\usepackage{psfrag}
\usepackage{eucal}
\usepackage{algorithm}
\usepackage{overpic}
\usepackage{subfigure}
\usepackage{url}
\usepackage[all,cmtip]{xy}
\usepackage{algorithm}
\usepackage[noend]{algpseudocode}

\definecolor{darkblue}{rgb}{0.0,0.0,0.6}
\hypersetup{colorlinks,breaklinks,
	linkcolor=darkblue,urlcolor=darkblue,
	anchorcolor=darkblue,citecolor=darkblue}
\newtheorem{assumption}{Assumption}
\newtheorem{definition}{Definition}
\newtheorem{lemma}{Lemma}
\newtheorem{corollary}{Corollary}
\newtheorem{proposition}{Proposition}
\newtheorem{theorem}{Theorem}
\newtheorem{remark}{Remark}

\makeatletter
\newcommand{\setassumptiontag}[1]{
  \let\oldtheassumption\theassumption
  \renewcommand{\theassumption}{#1}
  \g@addto@macro\endassumption{
    \addtocounter{assumption}{-1}
    \global\let\theassumption\oldtheassumption}
  }
\makeatother

\usepackage[normalem]{ulem} 

\newcommand{\pro}[2]{{\mathcal{P}_{#1}\left\{{#2}\right\}}}

\def \qed {\hfill \vrule height6pt width 6pt depth 0pt}

\newcommand{\R}{{\mathbb{R}}}
\newcommand{\bx}{{\mathbf{x}}}
\newcommand{\hbx}{\hat{\bx}}

\newcommand{\di}{{\mathrm{diag}}}

\newcommand{\dW}{{\|I-W\|^2}}

\newcommand{\dJi}[1]{{\nabla_{#1} J_{#1}}}

\newcommand{\A}{{\mathcal{A}}}

\newcommand{\Gra}{{\mathcal{G}}}

\newcommand{\tx}{{\tilde{x}}}

\newcommand{\Om}{{\Omega}}

\newcommand{\one}{{\mathbf{1}}}

\newcommand{\N}{{\mathcal{N}}}

\newcommand{\bF}{{\mathbf{F}}}

\newcommand{\trace}{{\mathrm{trace}}}
\newcommand{\T}{{\mathrm{T}}}

\def\an#1{{\color{black}#1}}
\def\BibTeX{{\rm B\kern-.05em{\sc i\kern-.025em b}\kern-.08em
		T\kern-.1667em\lower.7ex\hbox{E}\kern-.125emX}}

        \def\BibTeX{{\rm B\kern-.05em{\sc i\kern-.025em b}\kern-.08em
		T\kern-.1667em\lower.7ex\hbox{E}\kern-.125emX}}

\allowdisplaybreaks

\begin{document}

\begin{frontmatter}

\title{Fast Distributed Nash Equilibrium Seeking in Monotone Games} 

\thanks[footnoteinfo]{This paper was not presented at any IFAC 
meeting. Corresponding author T. Tatarenko. Tel. +49 6151 16-25048. 
}

\author[Paestum]{Tatiana Tatarenko}
\author[Rome]{Angelia Nedi\'c}

\address[Paestum]{TU Darmstadt, Germany}  
\address[Rome]{ASU, Arizona, USA}             


\begin{abstract}
This work proposes a novel distributed approach for \an{computing a Nash equilibrium} in convex games with merely monotone and restricted strongly monotone pseudo-gradients. By leveraging the idea of the centralized operator extrapolation method presented in \cite{LanExtrapolation} to solve variational inequalities, we develop the algorithm converging to Nash equilibria in games, where players have no access to the full information but are able to communicate  with neighbors over some communication graph. The convergence rate is demonstrated to be geometric and improves the rates obtained by the previously presented procedures seeking Nash equilibria in the class of games under consideration.

\end{abstract}

\end{frontmatter}

\section{INTRODUCTION}\label{sec:intro}

Game theory deals with a specific class of optimization problems arising in multi-agent systems, in which each
agent, also called player, aims to minimize its local cost function coupled through decision variables (actions) of all agents (players) in a system. 
The applications of game-theoretic optimization can be found, for example, in electricity markets, communication networks, autonomous driving systems and the smart grids \cite{Alpcan2005,BasharSG,GTM,Scutaricdma}. Solutions to such optimization problems are Nash equilibria  which characterize stable joint actions in games. To find these solutions in \an{a so called convex game, one can use their equivalent characterization as the solutions to the variational inequality} defined for the game's pseudo-gradient over the joint action set \cite{FaccPang1}. Moreover, it is known that, given a strongly monotone and Lipschitz continuous \an{mapping, the projection algorithm} converges geometrically fast to the unique solution of the variational inequality  and, thus, to the unique Nash equilibrium of the game. 
The convergence rate, in terms of the $k$th iterate's distance to the solution,
is in the order of $O\left(\exp\left\{-\frac{k}{\gamma^2}\right\}\right)$(see \cite{Nesterov}), where 
$\gamma = L/\mu\ge 1$ with $L$ and $\mu$ being the Lipschitz continuity and strong monotonicity constants of the mapping, respectively.
This rate has been improved in \cite{Nesterov} 
to the rate of $O\left(\exp\left\{-\frac{k}{\gamma}\right\}\right)$ 
by a more sophisticated algorithm
that requires, at each iteration, two operator evaluations and two projections. To relax these requirements, the paper \cite{LanExtrapolation} presents the so called operator extrapolation method achieving the same rate $O\left(\exp\left\{-\frac{k}{\gamma}\right\}\right)$ with one
operator evaluation and one projection per iteration. Moreover, geometrically fast convergence of the operator extrapolation method takes place under a weaker condition of restricted strong
monotonicity.  
These fast algorithms require full information in the sense that each player observes actions of all other players at every iteration.

Since in the modern large-scale
systems each agent has access only to some partial information about joint actions, \emph{fast distributed communication-based} optimization procedures in games have
gained a lot of attention over the recent years (see \cite{survey_distGT} for an extensive
review and bibliography). In particular, the work \cite{Bianchi2019} presents a proximal-point algorithm for converging to the Nash equilibrium with
a geometric rate in strongly monotone games. However, this algorithm requires the evaluation of a proximal operator, at each iteration, 
\an{which gives an implicit relation for the new iterate and the current one.} On the other hand, the papers \cite{ifac_TatNed} and \cite{directMethod_Grammatico} propose the distributed procedures based on the gradient algorithm and demonstrate their geometric convergence rate of the order $O\left(\exp\left\{-\frac{k}{\gamma^4}\right\}\right)$ for strongly monotone games with player communications
over time-invariant and time-varying graphs, respectively. The works \cite{Cdc2018_TatShiNed,AccGRANE_TAC} focus on a reformulation of a Nash equilibrium in distributed setting in terms of a so called augmented variational inequality, which takes into account the communication network that players are using. The main goal of such reformulations has been to adjust the fast centralized procedure from \cite{Nesterov} to the distributed settings and accelerate learning Nash equilibria under such settings. However, the acceleration (to the rate $O\left(\exp\left\{-\frac{k}{\gamma^3}\right\}\right)$) has been guaranteed only for a restrictive subclass of games with strongly monotone and Lipschitz continuous pseudo-gradients. 

It is worth noting that the distributed procedures mentioned above and their convergence rates have been investigated in the context of strongly monotone games. As for the case of merely monotone games, the seminal work~\cite{Nemir2005} presents the complexity result on the rate $O\left(\frac{1}{k}\right)$ that can be achieved by a so called \emph{gap function} in a mirror-prox extragradient method applied to smooth monotone variational inequality problems. Such method, however, involves iterates requiring the solution of two projection subproblems at each iteration. Moreover, it cannot be guaranteed that the mirror-prox extragradient method improves its performance when applied to solve smooth monotone variational inequalities. To overcome these limitations, the aforementioned work~\cite{LanExtrapolation} proposes a unified approach which not only estimates one projection operator per iterate but also demonstrates the optimal convergence rates in both merely monotone and (restricted) strongly monotone problems. 

Thus, inspired by the results achieved  for a \emph{centralized} solution of smooth monotone variational inequalities in~\cite{LanExtrapolation}, we present a novel fast \emph{distributed} discrete-time algorithm for seeking Nash equilibria in games with merely monotone and restricted strongly monotone pseudo-gradients. 
The developed procedure converges to the Nash equilibrium with the rate $O\left(\exp\left\{-\frac{k}{\gamma^2}\right\}\right)$  in the case of a restricted strongly monotone pseudo-gradient\footnote{To our knowledge, it is  the best known rate given the setting under consideration.}, whereas the gap function converges to zero with the rate $O\left(\frac{1}{k^{1/2-\epsilon}}\right)$, if the pseudo-gradient is merely monotone. To the best of our knowledge, convergence rate analysis of communication-based distributed procedures in merely monotone games has not been addressed so far.  

Some preliminary results of this paper have been presented in the conference version~\cite{ECC24_TatNed}. We emphasize that in the conference paper~\cite{ECC24_TatNed} the focus has been only on the restricted strongly monotone case and merely monotone games have not been studied. Moreover, due to the space limitation, the proof of the main technical result (Proposition~\ref{prop:3points}) has been omitted in~\cite{ECC24_TatNed}. In this paper, we provide the full proofs as well as the novel results on the convergence rate of the distributed communication-based procedure in merely monotone games. 

\textbf{Notations.}
The set $\{1,\ldots,n\}$ is denoted by $[n]$. Given a nonempty set $K\subseteq\R^n$ and a 
\an{mapping $f:K\to\R$, we let $\nabla_i f(x) = \frac{\partial f(x)}{\partial x_i}$ be} the partial derivative taken with respect to the $i$th coordinate of the vector variable $x\in\R^n$.
For any real vector space $\tilde E$ its dual space is denoted by $\tilde E^*$ and the inner product is denoted by $\langle u,v \rangle$, $u\in\tilde E^*$, $v\in \tilde E$.
A mapping $g:\tilde E\to \tilde E^*$ is said to be \emph{strongly monotone with the constant $\mu>0$ on the set $Q\subseteq \tilde E$}, if $\langle g(u)-g(v), u - v \rangle\ge\mu\|u - v\|^2$ for all $u,v\in Q$.
 A mapping $g:\tilde E\to \tilde E^*$ is said to be Lipschitz continuous on the set $Q\subseteq \tilde E$ with the constant $L$, if $\| g(u)-g(v)\|_*\le L\|u-v\|$ for all $u,v\in Q$, where \an{$\|\cdot\|_*$ is the norm in the dual space $\tilde E^*$.}
We consider real vector space $E$, which is either space of real vectors $E = E^* = \R^n$ or the space of real matrices $E = E^* = \R^{n\times n}$. 
In the case $E = \R^n$ we use $\|\cdot\|$ to denote the Euclidean norm induced by the standard dot product in $\R^n$. In the case $E = \R^{n\times n}$, 
the inner product $\langle u,v \rangle \triangleq \sqrt{\trace(u^Tv)}$ is the Frobenius inner product on $\R^{n\times n}$ and 
$\|\cdot\|$ denotes the Frobenius norm induced by the Frobenius inner product, i.e., $\|v\| \triangleq \sqrt{\trace(v^Tv)}$.
We use $\pro{\Om}{v}$ to denote the projection of $v\in E$ on a closed convex set $\Om\subseteq E$.
For any matrix $A$, the vector consisting of the diagonal entries of the matrix $A$ is denoted by $\di(A)$.

\section{Distributed Learning in Convex Games}

We consider a non-cooperative game between $n$ players. Let $J_i$ and $\Om_i\subseteq \R$ denote\footnote{All results below are applicable for games with different dimensions $\{d_i\}$ of the action sets $\{\Om_i\}$. The one-dimensional case is considered for the sake of notation simplicity.} respectively the cost function and the action set of the player $i$. We denote the joint action set by $\Om = \Om_1\times\cdots\times\Om_n$. Each function $J_i(x_i,x_{-i})$, $i\in[n]$, depends on $x_i$ and $x_{-i}$, where $x_i\in\Om_i$ is the action of the player $i$ and $x_{-i}\in\Om_{-i}=\Om_1\times\cdots\times\Om_{i-1}\times\Om_{i+1}\times\cdots\times\Om_n$ denotes the joint action of all players except for the player $i$. We assume that the players can interact over an undirected communication graph $\Gra([n],\A)$. The set of nodes is the set $[n]$ of players, and 
the set $\A$ of undirected arcs is such that $\{i,j\}\in\A$ whenever there is an undirected communication link between $i$ to $j$
and, thus, some information (message) can be passed between the players $i$ and $j$.
For each player $i$, the set $\N_i$ is the set of neighbors of player $i$ in the graph $\Gra([n],\A)$, i.e.,
$\N_{i}\triangleq\{j\in[n]: \, \{i,j\}\in\A\}$.
We denote such a game by $\Gamma(n,\{J_i\},\{\Om_i\},\Gra)$, and 
we make the following assumptions regarding the game.

\begin{assumption}\label{assum:convex}[Convex Game]
For all $i\in[n]$, the set $\Om_i$ is convex and closed, while the function $J_i(x_i, x_{-i})$ is convex and continuously differentiable in $x_i$ for each fixed $x_{-i}$.
\end{assumption}

When the cost functions $J_i(\cdot,x_{-i})$ are differentiable, we can define the pseudo-gradient.
\begin{definition}\label{def:gamemapping} The \emph{pseudo-gradient} $F(x):\Om\to\R^n$ of the game $\Gamma(n,\{J_i\},\{\Om_i\},\Gra)$ is defined as follows:
 $F(x)\triangleq\left[\nabla_1 J_1(x_1,x_{-1}), \ldots, \nabla_n J_n(x_n,x_{-n})\right]^T\in\R^{n}$,
 where $\nabla _i$ is the partial derivative with respect to $x_i$ (see \bf{Notations}).
\end{definition}

A solution to a game  is a Nash equilibrium, defined below.
\begin{definition}\label{def:NE}
 A vector $x^*=[x_1^*,x_2^*,\cdots, x_n^*]^T\in\Om$ is a \emph{Nash equilibrium} if for all $i\in[n]$ and all $x_i\in \Om_i$,
 $$J_i(x_i^*,x_{-i}^*)\le J_i(x_{i},x_{-i}^*).$$
 \end{definition}

By Assumption~\ref{assum:convex} and the connection between 
Nash equilibria and solutions to variational inequalities~\cite{FaccPang1}, 
the point $x^*\in\Om$ is a Nash equilibrium of the game $\Gamma(n,\{J_i\},\{\Om_i\},\Gra)$ if and only if the following
variational inequality holds
\begin{align}\label{eq:NE}
 \langle F(x^*), x-x^*\rangle\ge 0 \quad \mbox{for all $x\in\Om$}.
\end{align}
We make further assumptions regarding the players' action sets and cost functions, as follows.
\begin{assumption}\label{assum:compact}
The action sets $\Om_i$, $i\in[n]$, are bounded.
\end{assumption}
We notice that the existence of a Nash equilibrium is guaranteed under Assumptions~\ref{assum:convex} and~\ref{assum:compact} \cite{FaccPang1}. 
\begin{assumption}\label{assum:Lipschitz}
For every $i\in[n]$ the function $\dJi{i}(x_i,x_{-i})$ is Lipschitz continuous in
$x_i$ on $\Om_i$ for every fixed $x_{-i}\in\R^{n-1}$, that is, 
there exist a constant $L_i\ge 0$ such that  for all $x_{-i}\in\R^{n-1}$ we have for all $x_i,y_i\in\Om_i$,
\begin{align*}
 |\dJi{i}(x_i,x_{-i})-\dJi{i}(y_i,x_{-i})|&\leq L_i|x_i-y_i|.
\end{align*}
Moreover, for every $i\in[n]$ the function $\dJi{i}(x_i,x_{-i})$ is Lipschitz continuous in $x_{-i}$ on $\R^{n-1}$, for every fixed $x_i\in\Om_i$, that is, there is a constant $L_{-i}\ge 0$ such that 
for all $x_{i}\in\Om_i$ we have for all $x_{-i},y_{-i}\in \R^{n-1},$
\begin{align*}
|\dJi{i}(x_i,x_{-i})-\dJi{i}(x_i,y_{-i})|&\leq L_{-i}\|x_{-i}-y_{-i}\|.
\end{align*}
\end{assumption}
Notice that the assumption above is a standard one in the literature related to fast distributed Nash equilibrium seeking~\cite{AccGRANE_TAC,sergio19_prox,ifac_TatNed} and guarantees definite smooth properties of the problem allowing for development of a fast optimization procedure.  
\an{In particular, under Assumption~\ref{assum:Lipschitz}, the pseudo-gradient mapping $F(x)$ (see Definition~\ref{def:gamemapping}) is Lipschitz continuous on the set of vectors $x,y$ with $x_i,y_i\in\Om_i$ and $x_{-i},y_{-i}\in\mathbb{R}^{n-1}$ with the Lipschitz constant $L$ given by
\begin{equation}\label{eq-mapLip}
    L=\max_{i\in [n]}\sqrt{L_i^2+L_{-i}^2}.
    \end{equation}
This can be seen by using the analysis similar to the proof of Lemma~1 in~\cite{AccGRANE_TAC}.
}

The players' communications are restricted to the underlying connectivity graph $\Gra([n],\A)$, with which we associate a nonnegative symmetric mixing matrix $W$, i.e.,  a symmetric matrix with nonnegative entries and with positive entries $w_{ij}$ only when $\{i,j\}\in\A$. To ensure sufficient information  "mixing" in the network, we assume that the graph is connected. These assumptions are formalized, as follows.

\begin{assumption}\label{assum:connected}
The underlying undirected communication graph $\Gra([n],\A)$ is connected. The associated non-negative symmetric mixing matrix $W=[w_{ij}]\in\R^{n\times n}$ defines the weights on the undirected arcs such that $w_{ij}>0$  
if and only if $\{i,j\}\in\A$ and $\sum_{j=1}^{n}w_{ij} = 1$ for all $i\in[n]$.
\end{assumption}

\begin{remark}
There are some simple strategies for generating symmetric mixing matrices over undirected graphs for which Assumption~\ref{assum:connected} holds (see Section 2.4 in~\cite{Shi2014} for a summary of such strategies).
Recently, the work~\cite{ELLA} has addressed the problem of Nash equilibrium seeking in strongly monotone games endowed with time-varying directed communication. A geometric convergence rate has been obtained in terms of the weighted norms defined by the absolute probability sequence of the stochastic vectors relating the communication matrices, which prohibits an explicit computation of the convergence rate for the iterates themselves. 
\end{remark}
Assumption~\ref{assum:connected} implies that the second largest singular 
value $\sigma$ of $W$ is such that $\sigma\in(0, 1)$ and for any $x\in\R^n$ the following average property holds (see \cite{OlshTsits}):
\begin{align}\label{eq:sigma}
	\|Wx-\one\bar{x}\|\le \sigma\|x-\one\bar{x}\|,
\end{align}
where $\bar{x} = \frac{1}{n}\langle \one,x\rangle$ is the average of the elements of $x$. We will use the relation above in the analysis of the convergence rates (see the main technical result in Proposition~\ref{prop:3points}).  

In this work, we are interested in \emph{distributed seeking of the Nash equilibrium} in a game $\Gamma(n,\{J_i\},\{\Om_i\},\Gra)$ for which Assumptions~\ref{assum:convex}--\ref{assum:connected} hold. 
Moreover, we will analyze our approach under monotonicity condition on the pseudo-gradient. 
We will distinguish between two cases: the pseudo-gradient is \emph{restricted strongly monotone in respect to the Nash equilibrium $x^*$ over $\R^n$}, i.e.,
\begin{align}\label{eq:restr}
\langle F(x)-F(x^*), x - x^* \rangle\ge\mu\|x - x^*\|^2, \quad\forall x\in \R^n,
\end{align}
and the pseudo-gradient is \emph{monotone} over $\R^n$, i.e., 
\begin{align}\label{eq:mon}
\langle F(x)-F(y), x - y \rangle\ge0, \quad\forall x,y\in \R^n.
\end{align}
Note that the class of games with restricted strongly monotone pseudo-gradients  contains the class of games with strongly monotone pseudo-gradients. Moreover, given~\eqref{eq:restr} and existence of a Nash equilibrium, one concludes uniqueness of the solution (see, for example, Theorem 9 in~\cite{AccGRANE_TAC}).

In the next section, we present a distributed algorithm to find a Nash equilibrium in the game $\Gamma$ with the monotone pseudo-gradient.
\an{We will estimate the convergence rate of the algorithm in terms of the distance between the iterates and the solution when~\eqref{eq:restr} holds. For the case~\eqref{eq:mon}, we use 
the so called  \emph{gap function},  which is given by}
\begin{align}\label{eq:gap}
    g(y) = \max_{x\in\Omega} \langle F(x),y-x\rangle,
\end{align}
and \an{study its rate} of its convergence to $0$. 
\begin{remark}\label{rem:gap}
    The \an{gap function} above is continuous and satisfies $g(y) = 0$ if and only if $y$ is a Nash equilibrium in $\Gamma$, given that $F$ is continuous and Assumption~\ref{assum:convex} holds (see Theorem 1 in \cite{Nesterov}).

    In the main result on the convergence rate in merely monotone games, we will upper bound the gap function $g(y)$. Here we provide some intuition on what such a bound guarantees in terms of approximation to a Nash equilibrium. To do so, let us assume there is a point $y^*$ such that $g(y^*)<\epsilon$ for some $\epsilon>0$. According to the definition~\eqref{eq:gap}, the inequality
    $g(y^*)<\epsilon$ implies that 
    \begin{align}\label{eq:rem1}
    \langle F(x),x-y^*\rangle >-\epsilon, \quad \forall x\in\Om.        
    \end{align}
    \an{Assuming that $\epsilon\in(0,1)$, let us choose the point $x = \sqrt{\epsilon} z+(1-\sqrt{\epsilon} )y^*$, where $z\in\Om$ is arbitrary. By the convexity of $\Om$, we have $x\in\Om$ for any $z\in\Om$. Upon substituting this point} $x$ into~\eqref{eq:rem1}, we obtain
    \begin{align}\label{eq:rem2}
    \langle F(\sqrt{\epsilon} z+(1-\sqrt{\epsilon})y^*),z-y^*\rangle >-\sqrt{\epsilon},      
    \end{align}
    which implies that 
    \an{
    \begin{align}\label{eq:rem3}
        \langle F(y^*),z-y^*\rangle   >-(1+LD)\sqrt{\epsilon}, \quad \forall z\in\Om,
    \end{align}
    where $D = \max_{u,v\in\Om} \|v-u\|^2$ and $L$ is the Lipschitz continuity constant for the mapping $F$ due to Assumption~\ref{assum:Lipschitz} (see~\eqref{eq-mapLip}). To obtain the inequality in~\eqref{eq:rem3}, we used the following inequality to upper bound the left hand side of~\eqref{eq:rem2}: 
    \begin{align*}
        &\langle F(\sqrt{\epsilon} z+(1-\sqrt{\epsilon})y^*),z-y^*\rangle\cr
        &\le \langle F(y^*),z-y^*\rangle \cr
        &\qquad\qquad+ \langle F(\sqrt{\epsilon} z+(1-\sqrt{\epsilon})y^*) - F(y^*),z-y^*\rangle\cr
        &=\langle F(y^*),z-y^*\rangle+\sqrt{\epsilon} L\|z-y^*\|^2. 
    \end{align*}
    Next, using Assumption~\ref{assum:convex}, from~\eqref{eq:rem3} we obtain that for any $i\in[n]$ and any $z_i\in\Om_i$, 
    \[J_i(y_i^*,y_{-i}^*)\le J_i(z_i,y_{-i}^*) + (1+LD)\sqrt{\epsilon}.\]
    Thus, the condition~\eqref{eq:rem1} implies that $y^*$ is an $O(\sqrt{\epsilon})$-approximate Nash equilibrium.}   
\end{remark}  

\section{Algorithm Development}\label{sec:algo}

\subsection{Direct Acceleration}
Throughout the paper, we let player $i$ hold a \emph{local copy} of
the global decision\footnote{Note that global decision variable $x$ is a fictitious variable which never exists in the designed decentralized computing system.} variable $x$, which is denoted by
\begin{align}\label{eq:est}x_{(i)}=[\tx_{(i)1};\ldots;\tx_{(i)i-1};x_i;\tx_{(i)i+1};\ldots;\tx_{(i)n}]\in\R^n.
\end{align}
Here $\tx_{(i)j}$ can be viewed as a temporary estimate of $x_j$ by player $i$. With this notation, we always have $\tx_{(i)i}=x_i$. Also, we compactly denote the temporary estimates that player $i$ 
has for all decisions of the other players as
$$\tx_{-i}=[\tx_{(i)1};\ldots;\tx_{(i)i-1};\tx_{(i)i+1};\ldots;\tx_{(i)n}]\in\R^{n-1}.$$

We introduce the following estimation matrix:
$$
  \bx\triangleq\left(
     \begin{array}{ccc}
       & x_{(1)}^\T &  \\
       & x_{(2)}^\T &  \\
       &\vdots& \\
       & x_{(n)}^\T &  \\
     \end{array}
   \right)\in\R^{n\times n},
$$
where $x^\T$ denotes the transpose of a column-vector $x$.
We let $\Om_a$ denote an \an{augmented action set}, consisting of the estimation matrices, i.e., $\Om_a = \{\bx\in\R^{n\times n}\,|\, \di(\bx) \in \Om\}$.
The estimation of the pseudo-gradient of the game is defined as $\bF(\bx)\in\mathbb R^{n\times n}$, $\bx\in\Om_a$:
\begin{align}\label{eq:ext_gm}
   \bF(\bx)\!\triangleq\!\!\left(
     \begin{array}{cccc}
        \nabla_1 J_1(x_{(1)}) & 0 & \cdots &0\\
        0 & \nabla_2 J_2(x_{(2)}) & \cdots &0\\
        \vdots&\vdots&\ddots&\vdots\\
        0 & 0 & \cdots &\nabla_n J_n(x_{(n)})\\
     \end{array}
   \right).
\end{align}
The algorithm starts with an arbitrary initial point $\bx^0\in\Om_a$, that is, each player $i$ holds an arbitrary point $x_{(i)}^0\in\R^{i-1}\times\Om_i\times\R^{n-i}$. All the subsequent estimation matrices $\bx_{1}, \bx_{2},\ldots$ are obtained
through the updates described by Algorithm 1.

\begin{remark}
    Algorithm~1 is motivated by the operator extrapolation approach presented  in~\cite{LanExtrapolation} for solving variational inequalities in a centralized setting. The extrapolation here corresponds to the expression $\nabla_i J_i(\hat{x}_{(i)}^k) + \lambda_t[\nabla_i J_i({x}_{(i)}^k) - \nabla_i J_i(\hat{x}_{(i)}^{k-1})]$ in the update of the individual actions $x_i^{k+1}$, $i\in[n]$. It is inspired by the connection between Nesterov's acceleration and gradient extrapolation. In words, given a smooth optimization problem,  the gradient extrapolation approach can be viewed as a dual of the Nesterov’s accelerated gradient method (see~\cite{Lan2018} for more details). We refer to our proposed distributed algorithm \emph{Accelerated Direct Method} to emphasize that, in contrast to the work~\cite{AccGRANE_TAC}, \an{the communication step $\hat{x}_{(i)}^k = \sum\limits_{j\in\N_{i}} w_{ij} x_{(j)}^k$ \an{replaces} the corresponding centralized procedure that does not require an augmented game mapping.}
\end{remark}
\begin{center}
  {\textbf{Algorithm 1: Accelerated Direct Method}}
  \smallskip

    \begin{tabular}{l}
    \hline
   \emph{  } Set mixing matrix $W$;\\
    \emph{  } Choose step size $\alpha_t>0$ and parameter $\lambda_t>0$;\\
    \emph{  } Pick arbitrary $x^0_{(i)}\in\R^{i-1}\times\Om_i\times\R^{n-i}$, $i=1,\ldots,n$; \\
    \emph{  } Set $\hat{x}_{(i)}^0 = \sum\limits_{j\in\N_i} w_{ij} x_{(j)}^0$ and $x_{(i)}^1 = \hat{x}_{(i)}^0$,\\
    \emph{  } \textbf{for} $k=1,2,\ldots$, all players $i=1,\ldots,n$ do \\
    \emph{  } \quad$\hat{x}_{(i)}^k = \sum\limits_{j\in\N_i} w_{ij} x_{(j)}^k$, \\
    \emph{  } \quad$x_{i}^{k+1}=\mathcal{P}_{\Om_i}
    \left\{\hat x_i^k-\alpha_k\left[\nabla_i J_i(\hat{x}_{(i)}^k)\right.\right.$\\
    $\qquad\qquad\quad + \left.\left.\lambda_k[\nabla_i J_i({x}_{(i)}^k) - \nabla_i J_i(\hat{x}_{(i)}^{k-1})]\right]\right\}$;\\
    \emph{  } \qquad\textbf{for} $\ell=\{1,\ldots,i-1,i+1,\ldots,n\}$ \\
    \emph{  } \qquad\qquad$x_{(i)\ell}^{k+1}=\hat{x}_{(i)\ell}^k$;\\
    \emph{  } \qquad\textbf{end for};\\
    \emph{~ } \textbf{end for}.\\
    \hline
    \end{tabular}
\end{center}
The compact expression of Algorithm~1 in terms of the pseudo-gradient estimation $\bF(\cdot)$ (see~\eqref{eq:ext_gm}) and the estimation matrices $\{\bx^k\}$ is as follows: 
\begin{align}\label{eq:compactAlg}
\hbx^k &= W\bx^k,\\
\bx^{k+1} &= \pro{\Om_a}{\hbx^k - \alpha_k(\bF(\hbx^k)+\lambda_k(\bF(\bx^k) - \bF(\hbx^{k-1})))}.\qquad\nonumber
\end{align}
In the further analysis of Algorithm~1, we will use its compact matrix form as given in~\eqref{eq:compactAlg}.

\section{Analysis}
This subsection investigates  convergence  of Algorithm~1. To be able to prove its convergence to a unique solution and to state a linear rate, we first provide some important technical properties regarding the game pseudo-gradient $\bF$ and the iteration of the procedure.
\begin{lemma}\label{lem:Lip}
 Under Assumption~\ref{assum:Lipschitz}, the pseudo-gradient estimation $\bF(\cdot)$ is Lipschitz continuous over $\Om_a$ with the constant $L = \max_{i\in[n]}\sqrt{L^2_i + L^2_{-i}}$.
\end{lemma}
\begin{pf}
 See Lemma~1 in \cite{AccGRANE_TAC}. 
$\qed$\end{pf}

In the further analysis,  we will use the following lemma  (see Lemma 3.1 in \cite{LanBook} and the preceding discussion).
 \begin{lemma}\label{lem:3p}
     Let $Y\subseteq \R^{n\times n}$ be a closed convex set, and let 
     $y\in Y$ be defined through the following relation: 
     \[y:=\mathcal{P}_Y [z-\alpha_t g]\]
     for some $\alpha_t>0$ and  $g\in\R^{n\times n}$.   Then, we have 
     \begin{align*}
         \alpha_t\langle g,& y-x\rangle + \frac{1}{2}\|z-y\|^2\cr&
         \le \frac{1}{2}\|x-z\|^2 - \frac{1}{2}\|x-y\|^2\quad\hbox{for all $x\in Y$}.
    \end{align*}
 \end{lemma}


Before formulating the next result, let us define the following decomposition for each matrix $\bx\in\Om_a$: 
\begin{equation}\label{eq-perp-dec}
\bx = \bx_{||} + \bx_{\bot},\end{equation}
where $\bx_{||} = \frac{1}{n}\boldsymbol{1}\boldsymbol{1}^T\bx$ is the so called consensus matrix and $\bx_{\bot} = \bx - \bx_{||}$ with the implied property $\langle \bx_{||},\bx_{\bot} \rangle= 0$.
This decomposition will be used in the proof of the \an{forthcoming} results as it allows for an analysis of the gap function and for an estimation of the running weighted distance to the solution in the case of restricted strongly monotone pseudo-gradients. 
In the latter case, our aim is to upper bound a sum $\sum_{t=1}^k\omega_{1,t}\left\|\bx^{t+1}_{||}  -\bx\right\|^2 + \omega_{2,t}\|\bx^{t+1}_{\bot}\|^2$ by a sum $\sum_{t=1}^k\nu_{1,t}\left\|\bx^{t+1}_{||}  -\bx\right\|^2 + \nu_{2,t}\|\bx^{t+1}_{\bot}\|^2$ for some weights $\omega_{1,t}$, $\omega_{2,t}$, $\nu_{1,t}$, and $\nu_{2,t}$ in such a way that the final conclusion on linear convergence $\left\|\bx^{k+1}_{||}  -\bx\right\|^2 + \|\bx^{k+1}_{\bot}\|^2 =\left\|\bx^{k+1}  -\bx\right\|^2\le \rho^k\|\bx^{1}  -\bx\|$, $\rho\in(0,1)$ is implied.

The main technical result for the estimation of the convergence rates is provided in the next proposition. 
\begin{proposition}\label{prop:3points}
 Let Assumptions~\ref{assum:Lipschitz} and~\ref{assum:connected} hold. Let $\bx=\bx_{||}\in\Om_a$ be a consensus matrix with each row equal to some joint action $x$ in the game $\Gamma$. 
 Moreover, let positive sequences $\{\eta_t\}\in(0,1)$ and $\{\theta_t\}$ be such that $\theta_{t-1}\eta_{t}(1-\eta_{t-1})\ge \theta_tL^2\alpha_t^2\lambda_t^2$ for all $t\ge1$  and 
 $\theta_{t+1}\alpha_{t+1}\lambda_{t+1} = \theta_t\alpha_t$ for all $t\ge0$, 
 where $\lambda_t>0$ is the stepsize from Algorithm~1 and \an{$L = \max_{i\in[n]}\sqrt{L^2_i + L^2_{-i}}$}. Then, for any consensus matrix $\bx = \bx_{||}\in\Om_a$, we have
 \an{for all $k\ge1$,}
 \begin{align}\label{eq:3point}
 &\sum_{t=1}^{k} 
 \left[\theta_t \alpha_t\langle\bF(\bx^{t+1}),\bx^{t+1} -\bx\rangle+ \frac{\theta_t(1-\eta_t)}{2}\|\bx^{t+1} - \bx\|^2\right]\cr
 &\qquad\qquad-\frac{L^2\theta_k\alpha_k^2}{2(1-\eta_k)}\|\bx^{k+1}-\bx\|^2 \\
  &\le \sum_{t=1}^k [\theta_t b_{1,t} \left\|\bx^{t}_{||}  -\bx\right\|^2 + \theta_t b_{2,t}\|\bx^{t}_{\bot}\|^2]\nonumber, 
  \end{align}
with 
$b_{1,t} =\frac{1+ \eta_t\|W-I\|^2}{2}$ and $b_{2,t} = \frac{(1-\eta_t)\sigma^2 + \eta_t(1+\|W-I\|^2)}{2}$. 
\end{proposition}

\begin{pf}
Applying Lemma~\ref{lem:3p} to the iterates in \eqref{eq:compactAlg}, i.e., $y_t = \bx^t$, $Y = \Om_a$, and 
 $g_t= \bF(\hbx^t)+\lambda_t(\bF(\bx^t) - \bF(\hbx^{t-1}))$,
 we obtain for all $t\ge1$
 \begin{align*}
    & \alpha_t\left\langle \bF(\hbx^t)+\lambda_t(\bF(\bx^t) - \bF(\hbx^{t-1})), \bx^{t+1} - \bx\right\rangle \cr
     &\qquad+ \frac{1}{2}\|\bx^{t+1} - \hbx^t\|^2\le\frac{1}{2}\|\hbx^t-\bx\|^2  - \frac{1}{2}\|\bx^{t+1} - \bx\|^2,
 \end{align*} 
 {where $\bx\in\Om_a$ is any matrix from $\Om_a$ with the equal rows.}
  We multiply both sides of the preceding inequality by $\theta_t>0$, sum up the resulting relations over $t=1,\ldots, k,$ for an arbitrary $k\ge1$, to obtain
  {\allowdisplaybreaks
 \begin{align}\label{eq:3pl}
 &\sum_{t=1}^k[\theta_t\alpha_t\left\langle \bF(\hbx^t)+\lambda_t(\bF(\bx^t) - \bF(\hbx^{t-1})), \bx^{t+1} - \bx\right\rangle \cr
 &\qquad\qquad+ \frac{\theta_t}{2}\|\bx^{t+1} - \hbx^t\|^2]\cr
 &\le\sum_{t=1}^k\frac{\theta_t}{2}\left[\|\hbx^t-\bx\|^2  - \|\bx^{t+1} - \bx\|^2\right].
 \end{align}
  Next we consider the left hand side of the inequality above: 
 \begin{align}\label{eq:eq0}
  	&\sum_{t=1}^{k} \big[\theta_t \alpha_t\left\langle \bF(\hbx^t)+\lambda_t(\bF(\bx^t) - \bF(\hbx^{t-1})), \bx^{t+1} - \bx\right\rangle\cr
   &\qquad\qquad+ \frac{\theta_t}{2}\|\bx^{t+1} - \hbx^t\|^2\big]\cr
  	& = \sum_{t=1}^{k} 
	\big[\theta_t \alpha_t\langle \bF(\bx^{t+1}),\bx^{t+1} -\bx\rangle \cr
    &\quad -\theta_t\alpha_t \langle\bF(\bx^{t+1})-\bF(\hbx^{t}),\bx^{t+1} -\bx\rangle  \cr
    &\quad	+ \theta_t\alpha_t\lambda_t\langle\bF(\bx^{t})-\bF(\hbx^{t-1}),\bx^{t} -\bx\rangle \cr
  	&\quad+ \theta_t\alpha_t\lambda_t\langle\bF(\bx^{t})-\bF(\hbx^{t-1}),\bx^{t+1} -\bx^{t}\rangle + \frac{\theta_t}{2}\|\bx^{t+1} - \hbx^t\|^2\big]\cr
  	&=\sum_{t=1}^{k} 
	\big[\theta_t \alpha_t\langle\bF(\bx^{t+1}),\bx^{t+1} -\bx\rangle \cr
 &\quad+ \theta_t\alpha_t\lambda_t\langle\bF(\bx^{t})-\bF(\hbx^{t-1}),\bx^{t+1} -\bx^{t}\rangle + \frac{\theta_t}{2}\|\bx^{t+1} - \hbx^t\|^2\big]\cr
  	&\,-\theta_k\alpha_t \langle\bF(\bx^{k+1})-\bF(\hbx^{k}),\bx^{k+1} -\bx\rangle\cr
   &\,+\theta_0\alpha_t \langle\bF(\bx^{1})-\bF(\hbx^{0}),\bx^{1} -\bx\rangle\cr
   &=\sum_{t=1}^{k} 
	\big[\theta_t \alpha_t\langle\bF(\bx^{t+1}),\bx^{t+1} -\bx\rangle \cr
 &\quad+ \theta_t\alpha_t\lambda_t\langle\bF(\bx^{t})-\bF(\hbx^{t-1}),\bx^{t+1} -\bx^{t}\rangle + \frac{\theta_t}{2}\|\bx^{t+1} - \hbx^t\|^2\big]\\
   &\,-\theta_k\alpha_t \langle\bF(\bx^{k+1})-\bF(\hbx^{k}),\bx^{k+1} -\bx\rangle, \nonumber
 \end{align}
where in the second equality we used
 $\theta_{t+1}\alpha_{t+1}\lambda_{t+1} = \theta_t\alpha_t$ and in the last equality we used $\bx^{1}=\hbx^{0}$}. 
Next, we consider the sum of the second and the third term in relation~\eqref{eq:eq0}, for which we have
\begin{align}\label{eq:eq1}
 &\sum_{t=1}^{k} \!\left[\theta_t\alpha_t\lambda_t\langle\bF(\bx^{t})\!-\!\bF(\hbx^{t-1}),\bx^{t+1} \!-\!\bx^{t}\rangle + \frac{\theta_t}{2}\|\bx^{t+1} \!-\! \hbx^t\|^2\right]\cr
 &=\sum_{t=1}^{k} \left[\theta_t\alpha_t\lambda_t\langle\bF(\bx^{t})-\bF(\hbx^{t-1}),\bx^{t+1} -\bx^{t}\rangle \right.\cr
 &\quad\left.+ \frac{\theta_t\eta_t}{2}\|\bx^{t+1} - \hbx^t\|^2 +\frac{\theta_{t-1}(1-\eta_{t-1})}{2}\|\bx^{t} - \hbx^{t-1}\|^2\right]\cr
 &\qquad+ \frac{\theta_k(1-\eta_k)}{2}\|\bx^{k+1} - \hbx^k\|^2-\frac{\theta_{0}(1-\eta_0)}{2}\|\bx^{1} - \hbx^{0}\|^2\cr
 &=\sum_{t=1}^{k} \left[\theta_t\alpha_t\lambda_t\langle\bF(\bx^{t})-\bF(\hbx^{t-1}),\bx^{t+1} -\bx^{t}\rangle \right. \cr
 &\qquad\left.+\frac{\theta_t\eta_t}{2}\|\bx^{t+1} - \hbx^t\|^2+\frac{\theta_{t-1}(1-\eta_{t-1})}{2}\|\bx^{t} - \hbx^{t-1}\|^2\right]\cr
 &\qquad+ \frac{\theta_k(1-\eta_k)}{2}\|\bx^{k+1} - \hbx^k\|^2\cr
 & \ge \sum_{t=1}^{k} 
 \left[-L\theta_t\alpha_t\lambda_t\|\bx^{t}-\hbx^{t-1}\|\|\bx^{t+1} -\bx^{t}\|\right.\cr
 &\qquad\quad\left.+ \frac{\theta_t\eta_t}{2}\|\bx^{t+1} - \hbx^t\|^2 
 +\frac{\theta_{t-1}(1-\eta_{t-1})}{2}\|\bx^{t} - \hbx^{t-1}\|^2\right]\cr
 &\qquad + \sum_{t=1}^{k}\frac{\theta_t\eta_t}{2}
 \left[\|\bx^{t+1} - \hbx^t\|^2-\|\bx^{t+1} - \bx^t\|^2\right]\cr
 &\qquad+ \frac{\theta_k(1-\eta_k)}{2}\|\bx^{k+1} - \hbx^k\|^2\cr
 &\ge \sum_{t=1}^{k} 
 \left[-L\theta_t\alpha_t\lambda_t\|\bx^{t}-\hbx^{t-1}\|\|\bx^{t+1} -\bx^{t}\|\right. \cr
 &\qquad+ \left.\frac{\theta_t\eta_t}{2}\|\bx^{t+1} - \bx^t\|^2 +\frac{\theta_tL^2\alpha_t^2\lambda_t^2}{2\eta_t}\|\bx^{t} - \hbx^{t-1}\|^2\right]\cr
 &\qquad + \sum_{t=1}^{k}\frac{\theta_t\eta_t}{2}
 \left[\|\bx^{t+1} - \hbx^t\|^2-\|\bx^{t+1} - \bx^t\|^2\right]\cr
 &\qquad+ \frac{\theta_k(1-\eta_k)}{2}\|\bx^{k+1} - \hbx^k\|^2\cr
 &= \sum_{t=1}^{k} 
 (\|\bx^{t}-\hbx^{t-1}\|-\|\bx^{t+1} -\bx^{t}\|)^2\cr
 &\qquad + \sum_{t=1}^{k}\frac{\theta_t\eta_t}{2}
 \left[\|\bx^{t+1} - \hbx^t\|^2-\|\bx^{t+1} - \bx^t\|^2\right]\cr
 &\qquad+ \frac{\theta_k(1-\eta_k)}{2}\|\bx^{k+1} - \hbx^k\|^2\cr
 &\ge \sum_{t=1}^{k}\frac{\theta_t\eta_t}{2}[\|\bx^{t+1} - \hbx^t\|^2-\|\bx^{t+1} - \bx^t\|^2] \\
 &\qquad+ \frac{\theta_k(1-\eta_k)}{2}\|\bx^{k+1} - \hbx^k\|^2\nonumber,
\end{align}
where in the first equality we used the relation $\bx^{1} = \hbx^{0}$, and added and subtracted $(\theta_t\eta_t/2)\|\bx^{t+1} - \hbx^t\|^2$, in the first inequality we applied Lemma~\ref{lem:Lip} , whereas the second inequality is due to $\theta_{t-1}\eta_t(1-\eta_{t-1})\ge \theta_tL^2\alpha_t^2\lambda_t^2$. 
Next, we use again Lemma~\ref{lem:Lip}  to obtain
{\allowdisplaybreaks
\begin{align}\label{eq:last}
 &-\theta_k\alpha_k \langle\bF(\bx^{k+1})-\bF(\hbx^{k}),\bx^{k+1} -\bx\rangle \cr
 &\qquad+ \frac{\theta_k(1-\eta_k)}{2}\|\bx^{k+1} - \hbx^k\|^2\cr
 &\ge - \theta_k\alpha_k L \|\bx^{k+1}-\hbx^{k}\|\|\bx^{k+1} -\bx\| \cr
 &\qquad+ \frac{\theta_k(1-\eta_k)}{2}\|\bx^{k+1} - \hbx^k\|^2\cr
 &\ge -\frac{L^2\theta_k\alpha_k^2}{2(1-\eta_k)}\|\bx^{k+1}-\bx\|^2.
\end{align}}
The last inequality is obtained from $-2ab+a^2\ge-b^2$ with $a=\frac{\sqrt{\theta_k(1-\eta_k)}}{\sqrt{2}}\|\bx^{k+1} - \hbx^k\|$ and $b = \frac{L\sqrt{\theta_k}\alpha_k}{\sqrt{2(1-\eta_k)}}\|\bx^{k+1}-\bx\|$.
The preceding relations imply that
{\allowdisplaybreaks
\begin{align}\label{eq:eq2}
  	&\sum_{t=1}^{k} 
	\left[\theta_t \alpha_t\langle\bF(\bx^{t+1}),\bx^{t+1} -\bx\rangle + \frac{\theta_t}{2}\|\bx^{t+1} - \bx\|^2\right]\cr
 &\qquad-\frac{L^2\theta_k\alpha_k}{2(1-\eta_k)}\|\bx^{k+1}-\bx\|^2\cr
  &\stackrel{\eqref{eq:last}}{\le}\sum_{t=1}^{k} 
	\left[\theta_t \alpha_t\langle\bF(\bx^{t+1}),\bx^{t+1} -\bx\rangle + \frac{\theta_t}{2}\|\bx^{t+1} - \bx\|^2\right]\cr
 &\qquad-\theta_k\alpha_k \langle\bF(\bx^{k+1})-\bF(\hbx^{k}),\bx^{k+1} -\bx\rangle \cr
 &\qquad+ \frac{\theta_k(1-\eta_k)}{2}\|\bx^{k+1} - \hbx^k\|^2\cr
 &\stackrel{\eqref{eq:eq1}}{\le} \sum_{t=1}^{k} 
	\left[\theta_t \alpha_t\langle\bF(\bx^{t+1}),\bx^{t+1} -\bx\rangle + \frac{\theta_t}{2}\|\bx^{t+1} - \bx\|^2\right]\cr
 &\qquad-\theta_k\alpha_k \langle\bF(\bx^{k+1})-\bF(\hbx^{k}),\bx^{k+1} -\bx\rangle\cr
  &+ \sum_{t=1}^{k} \left[\theta_t\alpha_t\lambda_t\langle\bF(\bx^{t})-\bF(\hbx^{t-1}),\bx^{t+1} -\bx^{t}\rangle\right.\cr
 &\qquad\qquad\qquad\left. + \frac{\theta_t}{2}\|\bx^{t+1} - \hbx^t\|^2\right] \cr
  &\qquad - \sum_{t=1}^{k}\frac{\theta_t\eta_t}{2}[\|\bx^{t+1} - \hbx^t\|^2-\|\bx^{t+1} - \bx^t\|^2] \cr
  & \stackrel{\eqref{eq:eq0}}{=}\sum_{t=1}^{k} 
	\left[\frac{\theta_t}{2}\|\bx^{t+1} - \bx\|^2\right.\cr
 &\qquad\qquad\qquad\left.-\frac{\theta_t\eta_t}{2}[\|\bx^{t+1} - \hbx^t\|^2-\|\bx^{t+1} - \bx^t\|^2]\right]\cr
& \qquad + \sum_{t=1}^{k} \left[\theta_t \alpha_t\left\langle \bF(\hbx^t)+\lambda_t(\bF(\bx^t) - \bF(\hbx^{t-1})), \bx^{t+1} - \bx\right\rangle\right.\cr
&\qquad\qquad\qquad\left.+ \frac{\theta_t}{2}\|\bx^{t+1} - \hbx^t\|^2\right]\cr
   &\stackrel{\eqref{eq:3pl}}{\le} \sum_{t=1}^{k}
	\left[\frac{\theta_t}{2}\|\hbx^{t} - \bx\|^2 \right.\cr
 &\qquad\qquad\qquad\left.- \frac{\theta_t\eta_t}{2}\left[\|\bx^{t+1} - \hbx^t\|^2-\|\bx^{t+1} - \bx^t\|^2\right]\right]\cr
  	&=\sum_{t=1}^{k}
	\left[\frac{\theta_t(1-\eta_t)}{2}\|\hbx^{t} - \bx\|^2 \right.\\
 &\,\left.+ \frac{\theta_t\eta_t}{2}\left[\|\bx^{t+1} - \bx^t\|^2-\|\bx^{t+1} - \hbx^t\|^2+\|\hbx^{t} - \bx\|^2\right]\right].\nonumber
 \end{align}
}
Recalling our notation for $\bx_{||}$ and $\bx_\perp$ (see~\eqref{eq-perp-dec}), we have $\|\bx^{t+1}  -\bx\| \le \|\bx_{||}^{t+1}  -\bx\| + \|\bx_{\bot}^{t+1}\|$. 
As for the right hand side of~\eqref{eq:eq2}, we notice that 
\begin{align}\label{eq:rhs1}
 &\|\hbx^{t} - \bx\|^2 
 = \|W(\bx_{||}^{t} + \bx_{\bot}^{t}) - \bx\|^2 \cr
 &= \|\bx_{||}^{t}- \bx\|^2 + \|W\bx_{\bot}^t\|^2 \cr
 &\le \|\bx_{||}^{t}- \bx\|^2 + \sigma^2\|\bx_{\bot}^t\|^2,
\end{align}
where in the equality we used the fact that $\bx$ is a consensus matrix, whereas the inequality is due to~\eqref{eq:sigma}.
Moreover,
\begin{align}\label{eq:rhs2}
 &\|\bx^{t+1} - \bx^t\|^2-\|\bx^{t+1} - \hbx^t\|^2+\|\hbx^{t} - \bx\|^2 \cr
 &\le \|\bx^t - \bx\|^2 + \|\hbx^t-\bx^t\|^2 + \|\bx - \bx^{t+1}\|^2.
\end{align}
Combining~\eqref{eq:rhs1}-\eqref{eq:rhs2} with~\eqref{eq:eq2} and taking into account that 
\[\|\hbx^t-\bx^t\|^2 
= \|W\bx^t-\bx^t-W\bx+\bx\|^2 
= \|(I-W)(\bx^t-\bx)\|,\] 
we obtain 
\begin{align*}
 &\sum_{t=1}^{k} 
 \left[\theta_t \alpha_t\langle\bF(\bx^{t+1}),\bx^{t+1} -\bx\rangle+ \frac{\theta_t(1-\eta_t)}{2}\|\bx^{t+1} - \bx\|^2\right]\cr
 &\qquad\qquad\qquad-\frac{L^2\theta_k\alpha_k^2}{2(1-\eta_k)}\|\bx^{k+1}-\bx\|^2\cr
 &\le \sum_{t=1}^{k}\frac{\theta_t(1-\eta_t)}{2}
 \left[\|\bx_{||}^{t}- \bx\|^2 + \sigma^2\|\bx_{\bot}^t\|^2\right] \cr
 &\,+\sum_{t=1}^{k}\frac{\theta_t\eta_t}{2}
 \left[\|\bx_{||}^{t}- \bx\|^2 + \|\bx_{\bot}^t\|^2\right](1+\|(I-W)\|^2).
\end{align*}
After rearranging the terms, \an{we obtain the stated relation}. 
$\qed$\end{pf} 

\subsection{Monotone Case}
We proceed now with a refinement of the statement in Proposition~\ref{prop:3points} given monotonicity of the pseudo-gradient. 
\begin{lemma}\label{lem:3point_mon}
    Let the pseudo-gradient $F$ be monotone over $\R^n$ and the conditions of Proposition~\ref{prop:3points} hold. Then for any $x\in\Om$ the following holds: 
    \begin{align}\label{eq:3point_mon}
 &\sum_{t=1}^{k} 
 \left[\theta_t \alpha_t\langle F(x), x^{t+1} -x\rangle\right.\cr
 &\qquad\left. +\theta_tc_{1,t}\|\bx^{t+1}_{||}  -\bx \|^2+ \theta_tc_{2,t}\|\bx^{t+1}_{\bot}\|^2\right]\cr
 &\qquad\qquad\qquad-\frac{L^2\theta_k\alpha_k^2}{2(1-\eta_k)}\|\bx^{k+1}-\bx\|^2 \cr
  &\le \sum_{t=1}^k [\theta_t b_{1,t}  \|\bx^{t}_{||}  -\bx \|^2 + \theta_t b_{2,t}\|\bx^{t}_{\bot}\|^2], 
  \end{align}
  where $c_{1,t}=\frac{1-\eta_t}{2} - \alpha_t\zeta_t$ and $c_{2,t} = \frac{1-\eta_t}{2} - \alpha_t\left(L+\frac{L^2}{\zeta_t}\right)$, $b_{1,t} =\frac{1+ \eta_t\|W-I\|^2}{2}$, and $b_{2,t} = \frac{(1-\eta_t)\sigma^2 + \eta_t(1+\|W-I\|^2)}{2}$, with $\zeta_t$ being a sequence of positive numbers.
\end{lemma}

\begin{pf}
By adding and substracting $\bF(\bx^{t+1}_{||})$ to the right hand side of the inequality~\eqref{eq:3point} and using the relation $\bx^{t+1} = \bx_{||}^{t+1} +\bx_{\bot}^{t+1}$, we obtain 
\begin{align*}
&\langle\bF(\bx^{t+1}),\bx^{t+1} -\bx\rangle \cr
 &= \langle\bF(\bx^{t+1})-\bF(\bx^{t+1}_{||}) + \bF(\bx_{||}^{t+1}) ,
\bx_{||}^{t+1} +\bx_{\bot}^{t+1}  -\bx\rangle\cr
 &=\langle\bF(\bx^{t+1})- \bF(\bx_{||}^{t+1}),\bx^{t+1}  -\bx\rangle \cr
 &\qquad\qquad\qquad+ \langle
 \bF(\bx_{||}^{t+1})  ,\bx_{||}^{t+1}  -\bx\rangle + \langle\bF(\bx_{||}^{t+1})  ,\bx_{\bot}^{t+1}\rangle.\end{align*}
By Lemma~\ref{lem:Lip}, the mapping $\bF$ is Lipschitz continuous, implying that
\begin{align}\label{eq:Mon1}
&\langle\bF(\bx^{t+1}),\bx^{t+1} -\bx\rangle \ge-L\|\bx_{\bot}^{t+1}\|\|\bx^{t+1}  -\bx\|\cr
 &\qquad+\langle\bF(\bx_{||}^{t+1})  ,\bx_{||}^{t+1}  -\bx\rangle +  \langle\bF(\bx_{||}^{t+1})  ,\bx_{\bot}^{t+1}\rangle.
\end{align}
Next, we consider the term $\langle\bF(\bx_{||}^{t+1})  ,\bx_{||}^{t+1}  -\bx\rangle$.
\begin{align}\label{eq:Mon2}
&\langle\bF(\bx_{||}^{t+1})  ,\bx_{||}^{t+1}  -\bx\rangle = \langle F(x_{||}^{t+1}), x_{||}^{t+1}  -x\rangle\cr
&\ge \langle F(x), x_{||}^{t+1}  -x\rangle = \langle\bF(\bx)  ,\bx_{||}^{t+1}  -\bx\rangle\cr
&=\langle\bF(\bx)  ,\bx_{||}^{t+1} -\bx^{t+1}+ \bx^{t+1}-\bx\rangle\cr
& = - \langle\bF(\bx)  ,\bx_{\bot}^{t+1} \rangle+\langle F(x), x^{t+1}  -x\rangle,
\end{align}
where in the first inequality we used monotonicity of the pseudo-gradient $F$, whereas in the first and the last equalities we used the fact that $\langle\bF(\bx_{||}^{t+1})  ,\bx_{||}^{t+1}  -\bx\rangle = \langle F(x_{||}^{t+1}), x_{||}^{t+1}  -x\rangle$ and $\langle\bF(\bx),\bx^{t+1}-\bx\rangle = \langle F(x), x^{t+1}  -x\rangle$ respectively.  Thus, combining~\eqref{eq:Mon1} and~\eqref{eq:Mon2}, we obtain 
\begin{align*}
&\langle\bF(\bx^{t+1}),\bx^{t+1} -\bx\rangle  \ge-L\|\bx_{\bot}^{t+1}\|\|\bx^{t+1}  -\bx\|\cr &\qquad+\langle F(x), x^{t+1}  -x\rangle + \langle\bF(\bx_{||}^{t+1}) -\bF(\bx)  ,\bx_{\bot}^{t+1}\rangle\cr
&\ge -L\|\bx_{\bot}^{t+1}\|(\|\bx_{||}^{t+1}  -\bx\| + \|\bx_{\bot}^{t+1}\|)\cr 
&\qquad+\langle F(x), x^{t+1}  -x\rangle -L\|\bx_{||}^{t+1}  -\bx\|\|\bx_{\bot}^{t+1}\|\cr
& =\langle F(x), x^{t+1}  -x\rangle-2L\|\bx_{||}^{t+1}  -\bx\|\|\bx_{\bot}^{t+1}\| - L\|\bx_{\bot}^{t+1}\|^2\cr
&\ge\langle F(x), x^{t+1}  -x\rangle -\left(L+\frac{L^2}{\zeta_t}\right)\|\bx_{\bot}^{t+1}\|^2\cr
&\quad -\zeta_t\|\bx_{||}^{t+1}  -\bx\|^2,
\end{align*}
where in the last inequality we used the fact that $-2L\|\bx_{||}^{t+1}  -\bx\|\|\bx_{\bot}^{t+1}\|\ge -\frac{L^2}{\zeta_t}\|\bx_{\bot}^{t+1}\|^2
-\zeta_t\|\bx_{||}^{t+1}  -\bx\|^2$.
Finally, we apply this inequality to~\eqref{eq:3point} in Proposition~\ref{prop:3points} and use the fact that $\|\bx^{t+1}-\bx\|=\|\bx_{\bot}^{t+1}\|^2+\|\bx_{||}^{t+1}  -\bx\|^2$ to conclude the result.
$\qed$\end{pf}

Next, we formulate the main result for the case of merely monotone game. 
Let $\bar x^k$ be defined as the following weighted average of the joint actions obtained in the run of Algorithm~1: 
\begin{align}\label{eq:av}
    \bar x^k = \frac{\sum_{t=1}^k\theta_t\alpha_t x^{t+1}}{\sum_{t=1}^k\theta_t\alpha_t}.
\end{align}
Then the following statement holds for the sequence $\{\bar x^k\}$.
\begin{theorem}\label{th:main}
	Let the pseudo-gradient $F$ be monotone over $\R^n$, \an{Assumptions~\ref{assum:convex}--\ref{assum:connected} hold,} and the parameters in the Algorithm~1 be chosen as follows:
    $\alpha_k =\frac{A}{(k+1)^a}$, $\lambda_k = \left(\frac{k}{k+1}\right)^{a+b}$, where $a,b>0$, $a>1/2$, $a+b<1$, and $A< \frac{1}{2L}\frac{1}{2^{a+b/2}}$. Then, for the weighted average of the iterates $\bar x^k$ of the algorithm, see~\eqref{eq:av} with $\theta_t = \frac{1}{(t+1)^b}$,  the following holds: 
    \[\max_{x\in \Omega}\langle F(x), \bar x^k - x\rangle = O\left(\frac{1}{(k+1)^{1-a-b}}\right).\]
\end{theorem}
\begin{pf}
 First, we set up the parameters therein as follows: $\theta_t = \frac{1}{(t+1)^b}$, $\eta_t = \frac{1/2}{(t+1)^{2a-2\epsilon}}$, $\zeta_t = \frac{1}{(t+1)^{a-2\epsilon}}$, where $\epsilon$ is an arbitrary small positive number such that $2a-2\epsilon>1$. Note that conditions of Proposition~\ref{prop:3points}, namely $\theta_{t-1}\eta_{t}(1-\eta_{t-1})\ge \theta_tL^2\alpha_t^2\lambda_t^2$ for all $t\ge1$  and 
 $\theta_{t+1}\alpha_{t+1}\lambda_{t+1} = \theta_t\alpha_t$ for all $t\ge0$, hold for this choice of the parameters\footnote{Please, see Appendix~\ref{app1} for more details.}.
 Moreover, under this setting, there exists some finite $T\ge 1$ such that for any $t>T$ and $k>T$,
\begin{align}\label{eq:proof_mon1}
    &\theta_{t-1}c_{1,t-1}\ge\theta_tb_{1,t} \, \mbox{ and } \, \theta_{t-1}c_{2,t-1}\ge\theta_tb_{2,t}, \cr
    &\min\{c_{1,k}, c_{2,k}\} - \frac{L^2\alpha_k^2}{2(1-\eta_k)} \ge 0.
\end{align}
Please, refer to Appendix~\ref{app1} for more details on the inequalities above. Thus,~\eqref{eq:proof_mon1} and the iterate average $\bar x^k$ defined in~\eqref{eq:av}
combined with~\eqref{eq:3point_mon} imply  the following inequality: 
\begin{align}\label{eq:proof_mon2}
 &\left(\sum_{t=1}^{k} 
 \theta_t \alpha_t\right)\langle F(x), \bar x^{k} -x\rangle\cr
  &\le \sum_{t=1}^T [d_{1,t}  \|\bx^{t}_{||}  -\bx \|^2 + d_{2,t}\|\bx^{t}_{\bot}\|^2], 
  \end{align}
 where $d_{1,t} = \theta_tb_{1,t} - \theta_{t-1}c_{1,t-1}$ and $d_{2,t} = \theta_tb_{2,t} - \theta_{t-1}c_{2,t-1}$. 
 Finally, taking into account that $\Om$ is compact (Assumptions~\ref{assum:convex} and~\ref{assum:compact}) and $T$ is finite, we conclude that 
\begin{align*}
 &\max_{x\in\Omega}\langle F(x), \bar x^{k} -x\rangle = O\left(\frac{1}{(k+1)^{1-a-b}}\right), 
  \end{align*}
where we used the relation $\sum_{t=1}^{k} 
 \theta_t \alpha_t = O\left(\frac{1}{(k+1)^{1-a-b}}\right)$, which holds, given the settings for $\theta_t$ and $\alpha_t$.
$\qed$\end{pf}
The theorem above implies convergence of the gap function $g(\bar x^k)$ (see~\eqref{eq:gap}) to zero with a sublinear rate. Taking into account assumed compactness of $\Om$, continuity of $g$, and Remark~\ref{rem:gap}, we can conclude that any limit point of the sequence $\{\bar x^k\}$ is a Nash equilibrium in $\Gamma$. 

By optimizing the constants $a,b$ in Theorem~\ref{th:main}, we obtain the following result. 
\begin{corollary}\label{cor:results_mon}
	Let Assumptions~\ref{assum:convex}-\ref{assum:connected} hold, the pseudo-gradient $F$ be monotone over $\R^n$, and the parameters in Algorithm~1 be set up as follows: $\alpha_k = \frac{A}{(k+1)^{1/2+\epsilon/2}}$,  $A\le \frac{1}{2L}\frac{1}{\sqrt{2}}$, and $\lambda_k = \left(\frac{k}{k+1}\right)^{1/2+\epsilon}$, where $\epsilon$ is a positive number such that $\epsilon\in(0,1/2)$. Then there exists a Nash equilibrium in $\Gamma$ and any limit point of the sequence $\{\bar x^k\}$ defined by~\eqref{eq:av} with $\theta_t = \frac{1}{(t+1)^\epsilon}$ is a Nash equilibrium. Moreover, the following convergence rate for the gap function is obtained: 
    \[g(\bar x^k) = O\left(\frac{1}{k^{1/2-\epsilon}}\right).\] 
\end{corollary}
We notice that the result above provides the first rate estimation in distributed Nash equilibrium learning under the assumption of a merely monotone pseudo-gradient. In the light of Remark~\ref{rem:gap}, one can conclude that after $k$ iterates of Algorithm~1 the weighted average of the joint actions $\bar x^k$ achieves an $O\left(\frac{1}{t^{1/4-\epsilon}}\right)$-approximate Nash equilibrium in the game. Here we also refer to the paper~\cite{QuLiDist} for the result on accelerated distributed optimization, where the rate $O\left(\frac{1}{t^{1/4-\epsilon}}\right)$, $\epsilon\in(0,1/4)$, has been obtained for the case of convex and $L$-smooth objective functions. 
We would like to comment on the constant under big-$O$ notation in Corollary~\ref{cor:results_mon}. By choosing $\alpha_k = \frac{A}{(k+1)^{1/2+\epsilon/2}}$ and $\lambda_k = \left(\frac{k}{k+1}\right)^{1/2+\epsilon}$ with an arbitrary small $\epsilon$, we set up a small constant $b$ in the proof of Theorem~\ref{th:main}, which in its turn increases the finite $T$ (see~\eqref{eq:balance} in appendix related to the proof) and leads to a large constant under the big-$O$ notation (see~\eqref{eq:proof_mon2}).

The compact set $\Om$ \an{(Assumption~\ref{assum:convex} and Assumption~\ref{assum:compact})} allows us not only to conclude existence of a Nash equilibrium under the conditions of Theorem~\ref{th:main}, but also to take maximum over $\Om$ in~\eqref{eq:proof_mon2}. We notice that if Assumption~\ref{assum:compact} does not hold, but existence of a solution is guaranteed, one can analyze a so called restricted gap function \cite{Sedlmayer23a}.

\subsection{\an{Restricted} Strongly Monotone Case}
In this subsection, we focus on the case when the pseudo-gradient $F$ is restricted strongly monotone.
\begin{lemma}\label{lem:3point_str_mon}
    Let Assumptions~\ref{assum:convex},~\ref{assum:compact} hold and the pseudo-gradient $F$ be restricted strongly monotone over $\R^n$ in respect to the unique Nash equilibrium $x^*$ with some parameter $\mu>0$. Let the conditions of Proposition~\ref{prop:3points} hold. Then for a unique Nash equilibrium $x^*\in\Om$ the following holds: \an{for all $k\ge1$,}
    \begin{align}\label{eq:3point_str_mon}
 &\sum_{t=1}^k [\theta_t a_{1,t} \left\|\bx^{t+1}_{||}  -\bx^*\right\|^2 + \theta_t a_{2,t}\|\bx^{t+1}_{\bot}\|^2]  \cr
  &\qquad- \frac{\theta_kL^2\alpha_k^2}{2(1-\eta_k)}\|\bx^{k+1} - \bx^*\|^2 \cr
  &\le \sum_{t=1}^k [\theta_t b_{1,t} \left\|\bx^{t}_{||}  -\bx^*\right\|^2 + \theta_t b_{2,t}\|\bx^{t}_{\bot}\|^2], 
  \end{align}
where 
$a_{1,t} = \frac{1-\eta_t}{2} + \frac{\mu}{2n}\alpha_t$, $a_{2,t} = \frac{1-\eta_t}{2} - \left(L+\frac{2nL^2}{\mu}\right)\alpha_t$,	$b_{1,t} =\frac{1+ \eta_t\|W-I\|^2}{2}$, and $b_{2,t} = \frac{(1-\eta_t)\sigma^2 + \eta_t(1+\|W-I\|^2)}{2}$.
\end{lemma}

\begin{pf}
We focus on the right hand side of~\eqref{eq:3point} in Proposition~\ref{prop:3points}, where we set $\bx = \bx^*$. According to the definition~\eqref{eq:ext_gm} and taking into account that $x^*$ is the Nash equilibrium, we have $\langle\bF(\bx^*),\bx-\bx^*\rangle = \langle F(x^*), x-x^*\rangle\ge 0$ for any $\bx\in\Om_a$ (see~\eqref{eq:NE} and the definition of the mapping $\bF(\cdot)$ in~\eqref{eq:ext_gm}. 
Thus, we obtain
\begin{align*}
\langle\bF(\bx^{t+1}),\bx^{t+1} -\bx^*\rangle
&=\langle\bF(\bx^{t+1}) -\bF(\bx^*),\bx^{t+1} -\bx^*\rangle \cr
&\ +\langle\bF(\bx^*),\bx^{t+1}-\bx^*\rangle \cr
&\ge \langle\bF(\bx^{t+1})- \bF(\bx^*),\bx^{t+1} -\bx^*\rangle.
\end{align*}
According to our notation for $\bx_{||}$ and $\bx_\perp$ (see~\eqref{eq-perp-dec}), we have $\|\bx^{t+1}  -\bx^*\| \le \|\bx_{||}^{t+1}  -\bx^*\| + \|\bx_{\bot}^{t+1}\|$. 
Next, by adding and substracting $\bF(\bx^{t+1}_{||})$ to the right hand side of the preceding inequality and using the relation $\bx^{t+1} = \bx_{||}^{t+1} +\bx_{\bot}^{t+1}$, we obtain 
\begin{align*}
&\langle\bF(\bx^{t+1}),\bx^{t+1} -\bx^*\rangle \ge\cr
 & \langle\bF(\bx^{t+1})-\bF(\bx^{t+1}_{||}) + \bF(\bx_{||}^{t+1}) - \bF(\bx^*),\cr
&\qquad\qquad\qquad\qquad\qquad\qquad\qquad
\bx_{||}^{t+1} +\bx_{\bot}^{t+1}  -\bx^*\rangle\cr
 &=\langle\bF(\bx^{t+1})- \bF(\bx_{||}^{t+1}),\bx^{t+1}  -\bx^*\rangle \cr
 &\qquad+ \langle
 \bF(\bx_{||}^{t+1})  - \bF(\bx^*),\bx_{||}^{t+1}  -\bx^*\rangle\cr
 &\qquad + \langle\bF(\bx_{||}^{t+1})  - \bF(\bx^*),\bx_{\bot}^{t+1}\rangle.\end{align*}
By Lemma~\ref{lem:Lip}, the mapping $\bF$ is Lipschitz continuous, implying that
\begin{align}\label{eq:StrMon}
&\langle\bF(\bx^{t+1}),\bx^{t+1} -\bx^*\rangle \ge\cr
 &\ge -L\|\bx^{t+1}- \bx_{||}^{t+1}\|\|\bx^{t+1}  -\bx^*\|\cr
 &\qquad+\langle\bF(\bx_{||}^{t+1})  - \bF(\bx^*),\bx_{||}^{t+1}  -\bx^*\rangle \cr
 &\qquad- L\|\bx_{||}^{t+1}  - \bx^*\|\|\bx_{\bot}^{t+1}\|\cr
 & \ge -L\|\bx_{\bot}^{t+1}\|(\|\bx_{||}^{t+1}  -\bx^*\|+ \|\bx_{\bot}^{t+1}\|)\cr
 &\qquad + \langle \bF(\bx_{||}^{t+1})  - \bF(\bx^*),\bx_{||}^{t+1}  -\bx^*\rangle\cr
 &\qquad- L\|\bx_{||}^{t+1}  - \bx^*\|\|\bx_{\bot}^{t+1}\|\cr
 &\ge -2L\|\bx_{\bot}^{t+1}\|\|\bx_{||}^{t+1}  -\bx^*\| - L\|\bx_{\bot}^{t+1}\|^2 +\frac{\mu}{n}\|\bx_{||}^{t+1}  -\bx^*\|^2 \cr
 &\ge \frac{\mu}{2n}\|\bx_{||}^{t+1}  -\bx^*\|^2 - \left(L+\frac{2nL^2}{\mu}\right)\|\bx_{\bot}^{t+1}\|^2,
\end{align}
where in the last two inequalities we used Lemma~\ref{lem:Lip} and restricted strong monotonicity of $F$, implying that 
\begin{align*}
  &  \langle\bF(\bx_{||}^{t+1})  - \bF(\bx^*),\bx_{||}^{t+1}  -\bx^*\rangle \cr
&= \langle F(x_{||}^{t+1}) - F(x^*), x_{||}^{t+1}  -x^*\rangle\cr
&\ge \mu\|x_{||}^{t+1}  -x^*\|^2 
= \frac{\mu}{n}\|\bx_{||}^{t+1}  -\bx^*\|^2,
\end{align*}
and the fact that
\begin{align*}
    -2L&\|\bx_{\bot}^{t+1}\|\|\bx_{||}^{t+1}  -\bx^*\|\cr
    &\ge -\frac{2nL^2}{\mu}\|\bx_{\bot}^{t+1}\|^2 - \frac{\mu}{2n}\|\|\bx_{||}^{t+1}  -\bx^*\|^2.
\end{align*}
Using~\eqref{eq:StrMon} and rearranging the terms in~\eqref{eq:3point}, we conclude the result. 
$\qed$\end{pf}
We emphasize that compared to the work~\cite{LanExtrapolation}, where the centralized extrapolation-based procedure has been proposed, the above lemma requires restricted strong monotonicity over the whole $\R^n$ and not only over the joint action set $\Om$. This assumption is due to the fact that the players' estimates of joint actions must not belong to $\Om$ (see~\eqref{eq:est}). The analysis in the case of distributed settings, however, requires consideration of the estimation matrices.  

With Lemma~\ref{lem:3point_str_mon} at place we are ready to formulate the main result for games with restricted strongly monotone pseudo-gradients.
\begin{theorem}\label{th:main_str}
	Let Assumptions~\ref{assum:convex}-\ref{assum:connected} hold and the pseudo-gradient $F$ be restricted strongly monotone over $\R^n$ in respect to the unique Nash equilibrium $x^*$ with some parameter $\mu>0$. Let the parameters in the Algorithm~1 be chosen as follows: 
	\begin{align}\label{eq:g} 
		&\alpha_t = \alpha =  \min\{g_1, g_2, g_3, g_4\}, \, \lambda_t =\lambda = \frac{1}{1+\epsilon(\alpha_t)},
	\end{align}
where 
\begin{align*}
    &g_1 =  \frac{n\mu(1-\sigma^2)}{4(\mu+2nL)^2(1+\dW)},\cr
	&g_2 = \frac{n(1+\dW)}{2\mu},\cr
	&g_3 = \frac{\mu n(1+\dW)}{\mu^2+L^2(1+\dW)^2n^2}, \cr
	&g_4 = \frac{\mu}{\sqrt{4L^2\mu^2+16(L\mu+2nL^2)^2}},\cr
	&\epsilon(\alpha) = \frac{2\mu\alpha_t/n - (1+\dW)(1-\sqrt{1-4L^2\alpha_t^2})}{2+\dW(1-\sqrt{1-4L^2\alpha_t^2})}.	
\end{align*}
Then, $\epsilon(\alpha)>0$ and \an{for all $k\ge1$,}
\begin{align*}
   &\|\bx^{k+1} - \bx^*\|^2\cr
   &\le \frac{8+4\dW-4\dW\sqrt{1-4L^2\alpha^2}}{(1+\epsilon(\alpha))^{k-1}}\|\bx^1-\bx^*\|^2.
\end{align*}
\end{theorem}
\begin{pf}
  	Let $\theta_t=c^t$, where 
  	\begin{align}\label{eq:c}
  		c = \min\left\{\frac{a_1}{b_1},\frac{a_2}{b_2}\right\},
  	\end{align}
 and $a_1 = a_{1,t}$, $a_2=a_{2,t}$, $b_1=b_{1,t}$, and $b_2=b_{2,t}$ are defined in Proposition~\ref{lem:3point_str_mon}, namely,
 \begin{align*}
 	&a_1 = \frac{1-\eta}{2} + \frac{\mu}{2n}\alpha, a_2 = \frac{1-\eta}{2} - \left(L+\frac{2nL^2}{\mu}\right)\alpha, \cr
 	& b_1 =\frac{1+ \eta\|W-I\|^2}{2},\cr
  &b_2 = \frac{(1-\eta)\sigma^2 + \eta(1+\|W-I\|^2)}{2}. 
 \end{align*}
Moreover, let us choose $\lambda_t = \frac{\theta_t}{\theta_{t+1}}$ and $\eta = \frac{1-\sqrt{1-4L^2\alpha^2}}{2}$. 
Next, we demonstrate that under the condition~\eqref{eq:g}, $c>1$. For this purpose we check that in this case both $\frac{a_1}{b_1}$ and $\frac{a_2}{b_2}$ are larger than 1. 
Indeed, 
\begin{align*}
	\frac{a_1}{b_1} &= \frac{1-\eta+\mu\alpha/n}{1+\eta\|W-I\|^2},\cr
    \frac{a_2}{b_2} &= \frac{1-\eta-2\left(L+\frac{2nL^2}{\mu}\right)\alpha}{(1-\eta)\sigma^2 + \eta(1+\|W-I\|^2)}.
\end{align*}
First, we notice that under the condition $\alpha\le g_1$,  
$\frac{a_1}{b_1}\le\frac{a_2}{b_2}$ (see Appendix~\ref{app2}). Next, we have that 
$c=\frac{a_1}{b_1}>1$,
 if and only if 
 \begin{align}\label{eq:eta}
 	\eta&<\frac{\mu\alpha}{n(1+\|W-I\|^2)}.
 \end{align}
Given the definition of $\eta$, we conclude that 
under the condition $\alpha< \min\{g_2,g_3\}$,
\[\eta = \frac{1-\sqrt{1-4L^2\alpha^2}}{2} < \frac{\mu\alpha}{n(1+\|W-I\|^2)},\]
and, thus, \eqref{eq:eta} holds. 

Since $c>1$, it follows that
\begin{align*}
	L^2\alpha^2 = \eta(1-\eta)\le c \eta(1-\eta), 
\end{align*}
which implies $\theta_{t-1}\eta(1-\eta)\ge \theta_tL^2\alpha^2\lambda_t^2$. Hence, the conditions of Lemma~\ref{lem:3point_str_mon} hold and we conclude that 
 \begin{align*}
	&\theta_k[a_1 \|\bx_{||}^{k+1}  -\bx^*\|^2 +  a_2\|\bx_{\bot}^{t+1}\|^2]  - \frac{\theta_kL^2\alpha^2}{2(1-\eta)}\|\bx^{k+1} - \bx^*\|^2 \cr
 &\le \theta_1 [b_1 \|\bx_{||}^{1}  -\bx^*\|^2 +  b_2\|\bx_{\bot}^{1}\|^2], 
\end{align*}
where we used definition of $c$ in~\eqref{eq:c} implying $\theta_t b_1\le \theta_{t-1} a_1$ and $\theta_t b_2\le \theta_{t-1} a_2$. We also use the fact that $\bx_\bot=0$.
Thus, 
 $	\theta_k\min\{a_1,a_2\}\|\bx^{k+1}  -\bx\|^2  - \frac{\theta_kL^2\alpha^2}{2(1-\eta)}\|\bx^{k+1} - \bx\|^2 
	\le \theta_1 \max\{b_1,b_2\} \|\bx^{1}  -\bx\|^2.
$
As $\min\{a_1,a_2\} = a_2$ and $\max\{b_1,b_2\} = b_1$, 
we conclude that 
 \begin{align}\label{eq:main1}
	\theta_k\left[a_2- \frac{L^2\alpha^2}{2(1-\eta)}\right]\|\bx^{k+1}  -\bx^*\|^2  \le \theta_1 b_1 \|\bx^{1}  -\bx^*\|^2.
	\end{align}
Next, we notice that $a_2 = \frac{1-\eta}{2} - \left(L+\frac{2nL^2}{\mu}\right)\alpha$ is larger or equal to $1/8$, if the conditions $\alpha\le\frac{\mu}{4(L\mu+2nL^2)}$ and $\alpha\le\frac{\sqrt 3}{4L}$ hold (see Appendix~\ref{app2}).
On the other hand, given the condition $\alpha \le \frac{\sqrt 7}{8L}$, we have $\frac{L^2\alpha^2}{2(1-\eta)}\le \frac{1}{16}$. Taking this inequality together with $a_2\ge \frac{1}{8}$ into account, we obtain from~\eqref{eq:main1} that, given $\alpha\le\min\{g_4,g_5\}$, 
 \begin{align}\label{eq:main2}
	\frac{\theta_k}{16}\|\bx_{k+1}  -\bx^*\|^2  \le \theta_1 b_1 \|\bx_{1}  -\bx^*\|^2.
\end{align}
Finally, we use the fact that $\theta_k = c^k = (1+\epsilon(\alpha))^k$, where 
$\epsilon(\alpha) = \frac{a_1}{b_1} - 1 = \frac{2\mu\alpha/n - (1+\dW)(1-\sqrt{1-4L^2\alpha^2})}{2+\dW(1-\sqrt{1-4L^2\alpha^2})}$
to conclude the result from~\eqref{eq:main2}. 
 $\qed$\end{pf}
It is worth noting that, in contrast to Theorem~\ref{th:main} where the compact set $\Om$ allows for taking maximum over the joint action set in~\eqref{eq:proof_mon2}, the theorem above requires compactness of $\Om$ to merely guarantee the existence of a Nash equilibrium. This assumption can be, thus, replaced by any assumption implying existence of a Nash equilibrium (for example, strong monotonicity of the pseudo-gradient over $\Om$). 
By choosing optimal settings for the algorithm's parameter, we obtain the following result.
\begin{corollary}\label{cor:results_str_mon}
Let Assumptions~\ref{assum:convex}-\ref{assum:connected} hold and the pseudo-gradient $F$ be restricted strongly monotone over $\R^n$ in respect to the unique Nash equilibrium $x^*$ with some parameter $\mu>0$.
	Let the parameters be chosen as follows: $\alpha = O\left(\frac{\mu}{L^2n}\right)$, $\lambda_t = \lambda = O\left(1-\frac{\mu^2}{L^2n^2}\right)$.
    Then, the estimation matrix $\bx^k$ converges to the consensus matrix $\bx^*$ with the Nash equilibrium $x^*$ as its rows. Moreover, 
    \[\|\bx^k-\bx^*\|^2 = O\left(\exp\left\{-\frac{k}{\gamma^2n^2}\right\}\right).\] 
\end{corollary}   
\begin{remark}\label{rem:results_str_mon}
    We notice that the obtained convergence rate is faster in respect to both time $k$ and the game dimension $n$ than the rates of previously proposed methods for distributed learning of Nash equilibria in restricted strongly monotone games~\cite{AccGRANE_TAC,Bianchi2019,ifac_TatNed}. Indeed, the GRANE algorithm from~\cite{AccGRANE_TAC} is proven to converge to the Nash equilibrium with the rate $O\left(\exp\left\{-\frac{k}{\gamma^6n^6}\right\}\right)$, whereas the direct distributed procedure in~\cite{Bianchi2019,ifac_TatNed} improves this rate to $O\left(\exp\left\{-\frac{k}{\gamma^4n^3}\right\}\right)$.
\end{remark}

\section{Simulations}
\begin{figure}[!htb]
	\centering
	\includegraphics[width=0.5\textwidth]{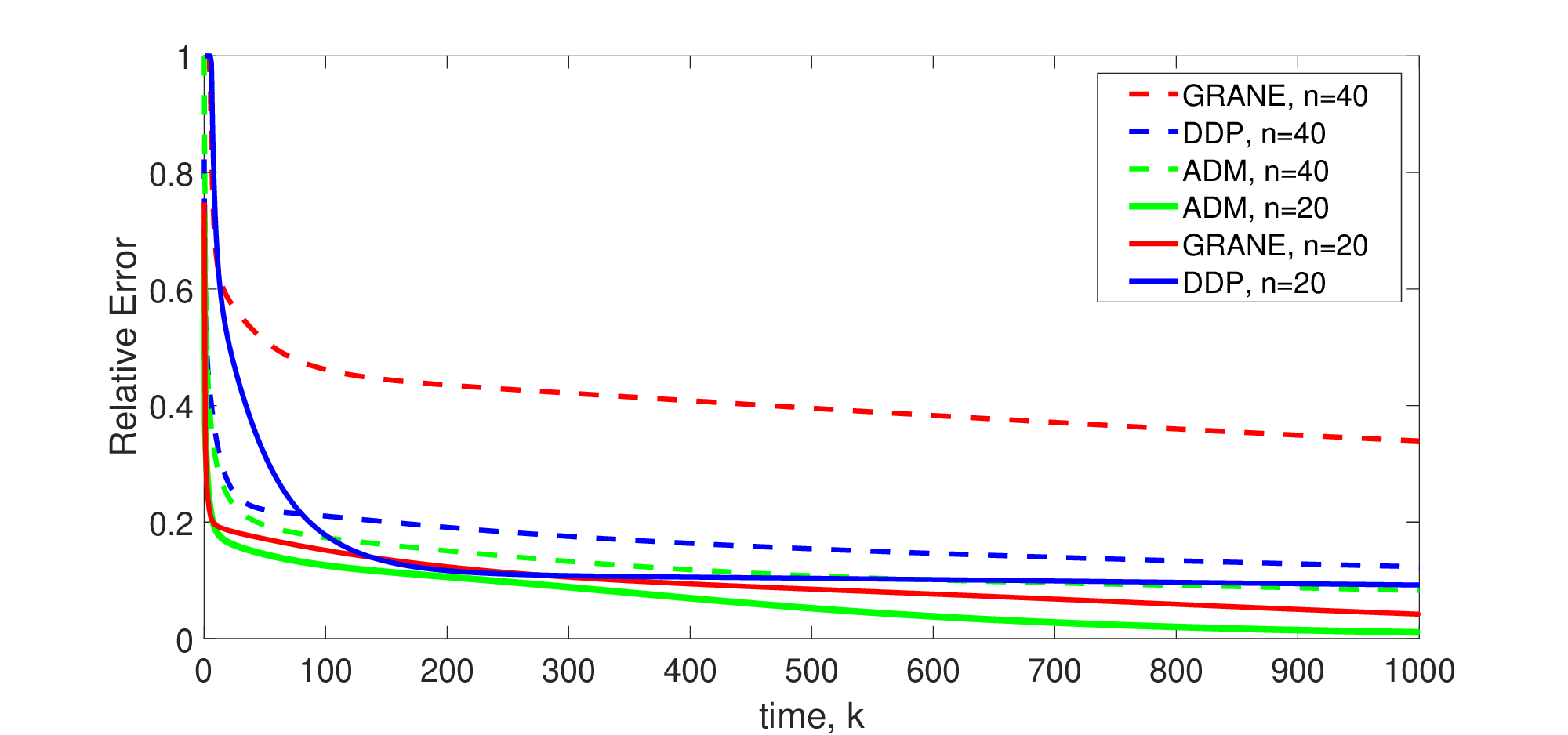}
	\caption{\label{eps:compare}Comparison of the presented accelerated direct method (ADM) with GRANE and DDP.}
\end{figure}
\begin{figure}[!htb]
	\centering
	\includegraphics[width=0.5\textwidth]{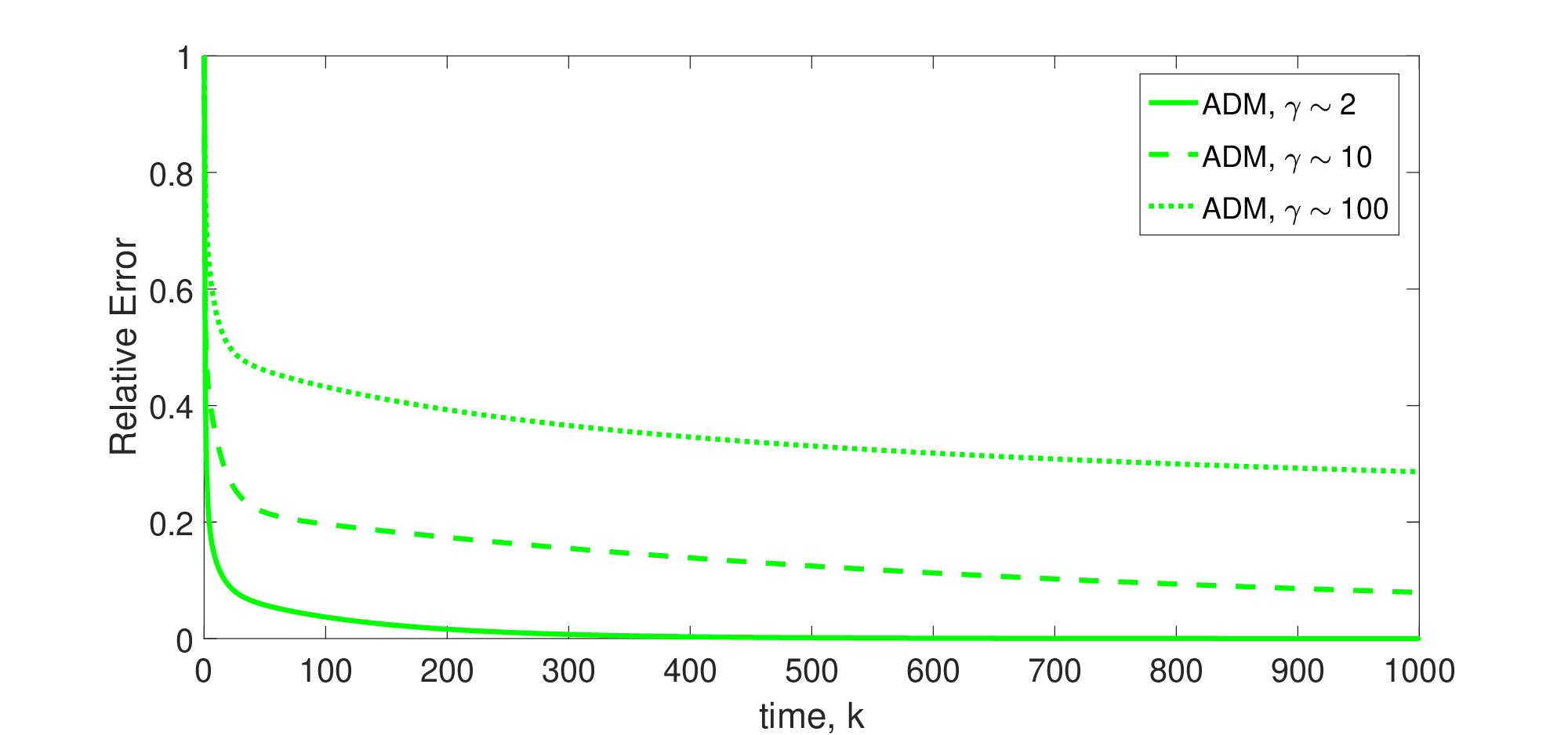}
	\caption{\label{fig:gamma}ADM in games with different condition numbers.}
\end{figure}
\begin{figure}[!htb]
	\centering
	\includegraphics[width=0.5\textwidth]{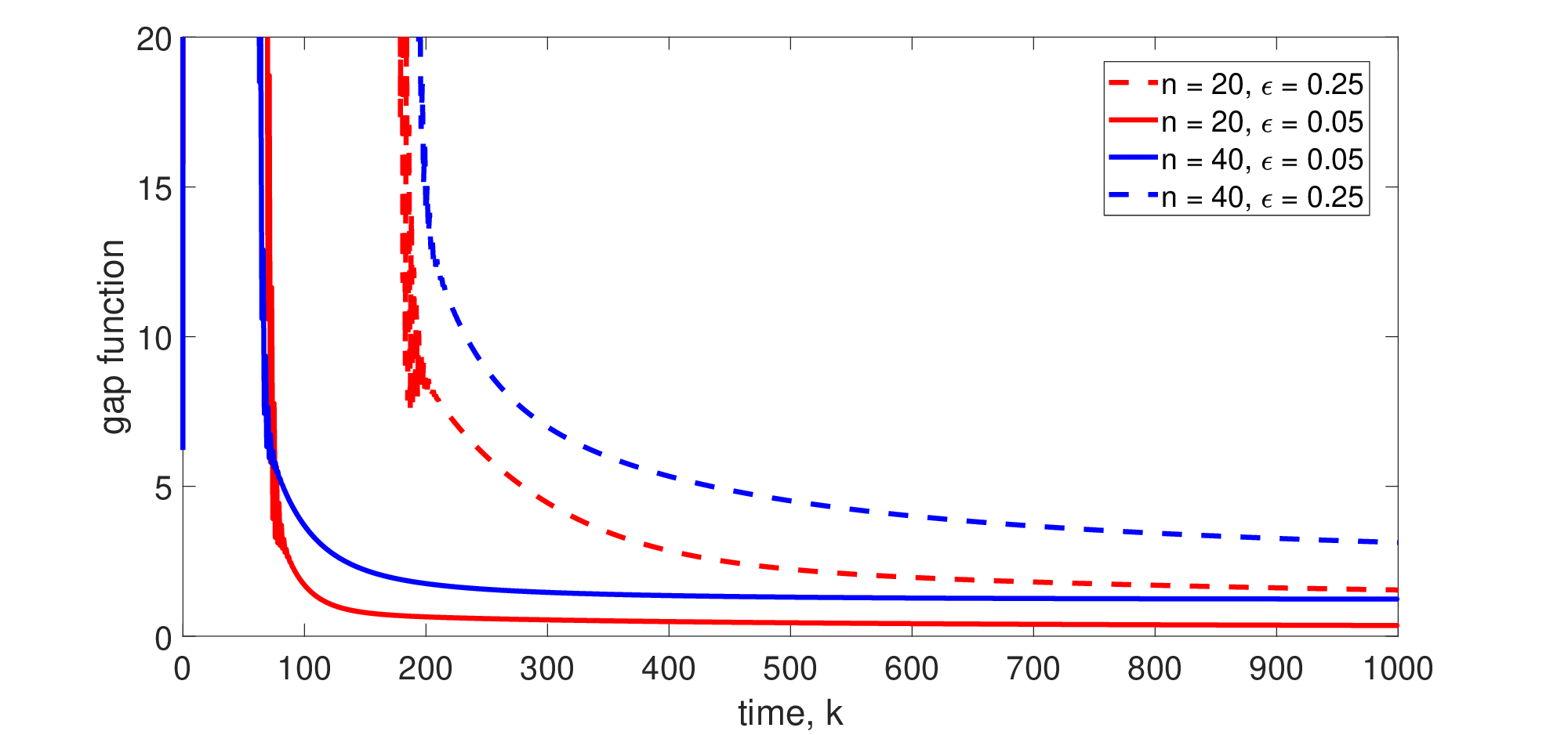}
	\caption{\label{fig:gap}ADM in merely monotone games.}
\end{figure}
Let us consider a class of games  with monotone pseudo-gradients. Specifically, we have players $\{1,2,\ldots,n\}$ and each player $i$'s objective is to minimize the cost function $J_i(x_i,x_{-i})=f_i(x_i)+l_i(x_{-i})x_i$, where $f_i(x_i)=0.5a_ix_i^2+b_ix_i$ and $l_i(x_{-i})=\sum_{j\neq i}c_{ij}x_j$. The local cost function is dependent on actions of all players, and the underlying communication graph is a randomly generated tree graph. 

First, we randomly select $a_i>0$, $b_i$, and $c_{ij}$ for all possible $i$ and $j$ to guarantee strong monotonicity of the pseudo-gradient.
We simulate the proposed gradient play algorithm and compare its implementation with the implementations of the algorithm GRANE presented in \cite{AccGRANE_TAC} and direct distributed procedure (DDP) from \cite{Bianchi2019,ifac_TatNed}. Figure~\ref{eps:compare} demonstrates the simulation results for games with $n=20,40$ (the dimension of the local action sets is equal to 2) and supports the theoretic conclusions formulated in Remark~\ref{rem:results_str_mon}. Moreover, Figure~\ref{fig:gamma} demonstrates the simulation result obtained by the ADM in  strongly monotone games with $n=20$ players under different condition numbers $\gamma = \frac{L}{\nu}$. As expected, the method slows down once $\gamma$ increases (see Corollary~\ref{cor:results_str_mon}). The evaluations are presented in terms of relative error equal to $\frac{\|x^k-x^*\|}{\|x^*\|}$, where the unique solution $x^*$ to the game is calculated by 20000 iterates of the centralized gradient play. 

Next, we simulate the ADM in the case of merely monotone games. To obtain this setting we randomly generate the parameters $a_i$, $b_i$ and $c_{ij}$, $i,j\in[n]$ and set up the matrix $C$ of the corresponding affine pseudo-gradient $F(x) = Cx + d$ to contain two equal rows. The final pseudo-gradient for implementation is $F_f(x) = C^TCx+b$.  Figure~\ref{fig:gap} provides the implementation result for the calculated gap function at the averaged iterates (see~\eqref{eq:av}), given the setting of Corollary~\ref{cor:results_mon} with $\epsilon = 0.05$ and $\epsilon=0.25$, $n=20$ and $n=40$ (the dimension of the local action sets is equal to 2). As we can see, the convergence rate slows down as $\epsilon$ increases. This is in consistence with Corollary~\ref{cor:results_mon}.  

\section{Conclusion}
This work extends centralized operator extrapolation method presented in \cite{LanExtrapolation} to distributed settings in games with \an{merely monotone and restricted strongly monotone pseudo-gradient mapping}, where players can exchange their information only with local neighbors via some communication graph.
In merely monotone games, a sublinear rate $O\left(\frac{1}{k^{1/2-\epsilon}}\right)$ is achieved for the value of the gap function. This is the first known result on the convergence rate analysis of distributed procedures applied to merely monotone games. In the case of restricted strongly monotone pseudo-gradient, the proposed procedure is proven to possess a geometric rate and to outperform the previously developed algorithms calculating Nash equilibria in games under the same assumptions. Future research directions include consideration of a more general communication topology and study of lower bounds for convergence rates of distributed methods in such class of games.

\bibliographystyle{plain}
\bibliography{accNE}

\begin{thebibliography}{10}

\bibitem{Alpcan2005}
T.~Alpcan and T.~Ba{\c{s}}ar.
\newblock {D}istributed {A}lgorithms for {N}ash {E}quilibria of {F}low
  {C}ontrol {G}ames.
\newblock In {\em Advances in Dynamic Games}, pages 473--498. Springer, 2005.

\bibitem{sergio19_prox}
M.~Bianchi, G.~Belgioioso, and S.~Grammatico.
\newblock A distributed proximal-point algorithm for nash equilibrium seeking
  under partial-decision information with geometric convergence.
\newblock {\em arXiv preprint arXiv:1910.11613}, 2019.

\bibitem{Bianchi2019}
M.~Bianchi, G.~Belgioioso, and S.~Grammatico.
\newblock A fully-distributed proximal-point algorithm for nash equilibrium
  seeking with linear convergence rate.
\newblock In {\em 2020 59th IEEE Conference on Decision and Control (CDC)},
  pages 2303--2308, 2020.

\bibitem{directMethod_Grammatico}
M.~Bianchi and S.~Grammatico.
\newblock Fully distributed nash equilibrium seeking over time-varying
  communication networks with linear convergence rate.
\newblock {\em IEEE Control Systems Letters}, 5(2):499--504, 2021.

\bibitem{LanExtrapolation}
G.~Kotsalis, G.~Lan, and T.~Li.
\newblock Simple and optimal methods for stochastic variational inequalities,
  i: Operator extrapolation.
\newblock {\em SIAM Journal on Optimization}, 32(3):2041--2073, 2022.

\bibitem{LanBook}
G.~Lan.
\newblock First-order and stochastic optimization methods for machine learning.
\newblock Springer, 2020.

\bibitem{Lan2018}
G.~Lan and Y.~Zhou.
\newblock Random gradient extrapolation for distributed and stochastic
  optimization.
\newblock {\em SIAM Journal on Optimization}, 28(4):2753--2782, 2018.

\bibitem{GTM}
N.~Li, Y.~Yao, I.~Kolmanovsky, E.~Atkins, and A.~R. Girard.
\newblock Game-theoretic modeling of multi-vehicle interactions at uncontrolled
  intersections.
\newblock {\em IEEE Transactions on Intelligent Transportation Systems},
  23(2):1428--1442, 2022.

\bibitem{Nemir2005}
A.~Nemirovski.
\newblock Prox-method with rate of convergence o(1/t) for variational
  inequalities with lipschitz continuous monotone operators and smooth
  convex-concave saddle point problems.
\newblock {\em SIAM Journal on Optimization}, 15(1):229--251, 2004.

\bibitem{Nesterov}
Yu. Nesterov and L.~Scrimali.
\newblock Solving strongly monotone variational and quasi-variational
  inequalities.
\newblock {\em Discrete and Continuous Dynamical Systems - A},
  31(4):1383--1396, 2011.

\bibitem{ELLA}
D.T.A. Nguyen, D.T. Nguyen, and A.~Nedi\'c.
\newblock Distributed nash equilibrium seeking over time-varying directed
  communication networks.
\newblock {\em IEEE Transactions on Control of Network Systems}, pages 1--12,
  2025.

\bibitem{OlshTsits}
A.~Olshevsky and J.~Tsitsiklis.
\newblock Convergence speed in distributed consensus and averaging.
\newblock {\em SIAM Journal on Control and Optimization}, 48(1):33--55, 2009.

\bibitem{FaccPang1}
J.-S. Pang and F.~Facchinei.
\newblock {\em Finite-dimensional variational inequalities and complementarity
  problems : vol. 1}.
\newblock Springer series in operations research. Springer, New York, Berlin,
  Heidelberg, 2003.

\bibitem{QuLiDist}
G.~Qu and N.~Li.
\newblock Accelerated distributed {N}esterov gradient descent.
\newblock {\em IEEE Transactions on Automatic Control}, 65(6):2566--2581, 2020.

\bibitem{BasharSG}
W.~Saad, H.~Zhu, H.~V. Poor, and T.~Ba{\c{s}}ar.
\newblock Game-theoretic methods for the smart grid: An overview of microgrid
  systems, demand-side management, and smart grid communications.
\newblock {\em IEEE Signal Processing Magazine}, 29(5):86--105, 2012.

\bibitem{Scutaricdma}
G.~Scutari, S.~Barbarossa, and D.~P. Palomar.
\newblock Potential games: A framework for vector power control problems with
  coupled constraints.
\newblock In {\em 2006 IEEE International Conference on Acoustics Speech and
  Signal Processing Proceedings}, volume~4, pages 241--244, May 2006.

\bibitem{Sedlmayer23a}
M.~Sedlmayer, D.-K. Nguyen, and R.~I. Bot.
\newblock A fast optimistic method for monotone variational inequalities.
\newblock In Andreas Krause, Emma Brunskill, Kyunghyun Cho, Barbara Engelhardt,
  Sivan Sabato, and Jonathan Scarlett, editors, {\em Proceedings of the 40th
  International Conference on Machine Learning}, volume 202 of {\em Proceedings
  of Machine Learning Research}, pages 30406--30438. PMLR, 23--29 Jul 2023.

\bibitem{Shi2014}
W.~Shi, Q.~Ling, G.~Wu, and W.~Yin.
\newblock {EXTRA}: {A}n {E}xact {F}irst-{O}rder {A}lgorithm for {D}ecentralized
  {C}onsensus {O}ptimization.
\newblock {\em SIAM Journal on Optimization}, 25(2):944--966, 2015.

\bibitem{ifac_TatNed}
T.~Tatarenko and A.~Nedić.
\newblock Geometric convergence of distributed gradient play in games with
  unconstrained action sets.
\newblock {\em IFAC-PapersOnLine}, 53(2):3367--3372, 2020.
\newblock 21st IFAC World Congress.

\bibitem{ECC24_TatNed}
T.~Tatarenko and A.~Nedi\'c.
\newblock Accelerating distributed nash equilibrium seeking.
\newblock In {\em 2024 European Control Conference (ECC)}, pages 323--328,
  2024.

\bibitem{Cdc2018_TatShiNed}
T.~{Tatarenko}, W.~{Shi}, and A.~{Nedić}.
\newblock Accelerated gradient play algorithm for distributed nash equilibrium
  seeking.
\newblock In {\em 2018 IEEE Conference on Decision and Control (CDC)}, pages
  3561--3566, 2018.

\bibitem{AccGRANE_TAC}
T.~Tatarenko, W.~Shi, and A.~Nedi\'c.
\newblock Geometric convergence of gradient play algorithms for distributed
  nash equilibrium seeking.
\newblock {\em IEEE Transactions on Automatic Control}, 66(11):5342--5353,
  2021.

\bibitem{survey_distGT}
M.~Ye, Q.-L. Han, L.~Ding, and S.~Xu.
\newblock Distributed nash equilibrium seeking in games with partial decision
  information: A survey.
\newblock {\em Proceedings of the IEEE}, 111(2):140--157, 2023.

\end{thebibliography}

\appendix
\section{More details on proof of Theorem~\ref{th:main}.}\label{app1}
Since $\alpha_t = \frac{A}{(t+1)^a}$, $\lambda_t = \left(\frac{t}{t+1}\right)^{a+b}$, $\theta_t = \frac{1}{(t+1)^b}$, we get 
$\theta_{t+1}\alpha_{t+1}\lambda_{t+1} = \theta_t\alpha_t$ for all \an{$t\ge 0$}.
Moreover, as $\eta_t = \frac{1/2}{(t+1)^{2a-\epsilon}}$, the inequality $\theta_{t-1}\eta_{t}(1-\eta_{t-1})\ge \theta_tL^2\alpha_t^2\lambda_t^2$ holds for all $t\ge1$. Indeed, given that $\lambda_{t} = \frac{\theta_{t-1}\alpha_{t-1}}{\theta_{t}\alpha_{t}}$, we obtain
  \begin{align*}
  &\frac{\theta_{t-1}\eta_{t}(1-\eta_{t-1})}{\theta_tL^2\alpha_t^2\lambda_t^2} =  \frac{\theta_{t}\eta_{t}(1-\eta_{t-1})}{\theta_{t-1}L^2\alpha_{t-1}^2}  \cr
  &=\left(\frac{t}{t+1}\right)^bt^{2\epsilon} \frac{2t^{2a-2\epsilon}-1}{4L^2A^2(t+1)^{2a-2\epsilon}}\ge 1,
  \end{align*}
  where the last inequality follows from the fact that
  \[\left(\frac{t}{t+1}\right)^bt^{2\epsilon}\ge 4L^2A^2\frac{(t+1)^{2a-2\epsilon}}{2t^{2a-2\epsilon}-1},\]
  as $\left(\frac{t}{t+1}\right)^bt^{2\epsilon}$ is an increasing function in \an{$t$ for $t\ge 1$}, whereas $\frac{(t+1)^{2a-2\epsilon}}{2t^{2a-2\epsilon}-1}\ge \frac{(t+2)^{2a-2\epsilon}}{2(t+1)^{2a-2\epsilon}-1}$ for $t=1,2,\ldots,$ and $A\le \frac{1}{2L}\frac{1}{2^{a+b/2-\epsilon}}$.

To obtain~\eqref{eq:proof_mon1}, we notice that under the choice of $\zeta_t = \frac{1}{(t+1)^{a-2\epsilon}}$, one \an{can conclude} that 
\[\frac{c_{1,t-1}}{b_{1,t}}\sim 1 - O\left(\frac{1}{t^{2a-2\epsilon}}\right),\]
whereas 
\begin{align}\label{eq:balance}
    \frac{\theta_t}{\theta_{t-1}}\sim 1 - \frac{b}{t}.
\end{align}
Since $2a-2\epsilon>1$, there exists $T_1$ such that $\theta_{t-1}c_{1,t-1}\ge\theta_tb_{1,t}$ for any $t\ge T_1$. The inequality $\theta_{t-1}c_{1,t-1}\ge\theta_tb_{1,t}$ holds for any sufficiently large $t\ge T_2$ ($T_2>0$) as $\sigma\in(0,1)$. 
Finally, $\min\{c_{1,k},c_{2,k}\} = c_{2,k} = \frac{1-\eta_k}{2} - \alpha_k\left(L+\frac{L^2}{\zeta_k}\right)$, which converges to $1/2$ as $k\to\infty$. Thus, there exists $T_3$ such that  $\min\{c_{1,k},c_{2,k}\}\ge \frac{L^2\alpha_k^2}{2(1-\eta_k)}$ for $k\ge T_3$, as $\frac{L^2\alpha_k^2}{2(1-\eta_k)}$ decreases monotonically to 0 as $k\to\infty$. By taking $T=\max\{T_1,T_2,T_3\}$ we obtain~\eqref{eq:proof_mon1}.

\section{More details on proof of Theorem~\ref{th:main_str}.}\label{app2}
\emph{The inequality $\frac{a_1}{b_1}\le\frac{a_2}{b_2}$ under the condition $\alpha_t<g_1$.}
{\allowdisplaybreaks
The condition $\alpha_t<g_1 = \frac{n\mu(1-\sigma^2)}{4(\mu+2nL)^2(1+\dW)}$ implies
$
 \alpha_t	\left(1+\frac{2nL}{\mu}\right)^2 \le \frac{n(1-\sigma^2)}{4\mu(1+\dW)}.
$
Thus, since $\eta_t = \frac{1-\sqrt{1-4L^2\alpha_t^2}}{2}<\frac{1}{2}$, we conclude that $(1-\eta_t)^2>\frac{1}{4}$. Hence,
\begin{align*}
	\alpha_t	\left(1+\frac{2nL}{\mu}\right)^2(1+\dW) \le \frac{n}{\mu}(1-\sigma^2)(1-\eta_t)^2.
\end{align*}
Next, we use $\eta_t<1$ and $\sigma<1$ to obtain
\begin{align*}&\left(1+\frac{2nL}{\mu}\right)^2(1+\dW) \cr
&\ge \left(1-\frac{2nL}{\mu}\right)^2(1+\eta_t\dW) - (1-\eta_t)(1-\sigma^2).
\end{align*}
Combining two last inequalities, we get 
\begin{align*}
	&\alpha_t	\left[\left(1+\frac{2nL}{\mu}\right)^2(1+\eta_t\dW)- (1-\eta_t)(1-\sigma^2)\right] \cr
 &\le \frac{n}{\mu}(1-\sigma^2)(1-\eta_t)^2.
\end{align*}
By multiplying both sides by $\frac{\mu}{n}$, we obtain
\begin{align*}
	&\alpha_t \left[ \left(\frac{\mu}{n}+2L+\frac{4nL^2}{\mu}\right)(1+\eta_t\dW)\right.\cr
 &\qquad\left.- (1-\eta_t)(1-\sigma^2)\frac{\mu}{n}\right] \le (1-\sigma^2)(1-\eta_t)^2.
 \end{align*}
 By adding $(1-\eta_t)(1+\eta_t\dW)$ to both sides of the preceding inequality and rearranging the terms, we have
    \begin{align*}
	&(1-\eta_t)(1+\eta_t\dW)+\alpha_t\frac{\mu}{n}(1+\eta_t\dW) \cr
 &\qquad\qquad- \frac{\mu}{n}\alpha_t(1-\eta_t)(1-\sigma^2)- (1-\eta_t)^2(1-\sigma^2)\cr
	&\le
	(1-\eta_t)(1+\eta_t\dW)\cr
 &\qquad\qquad-2\left(L+\frac{2nL^2}{\mu}\right)\alpha_t(1+\eta_t\dW).\end{align*}
 By further grouping the common terms, we obtain
    \begin{align*}
	&(1+\eta_t\dW)\left(1-\eta_t+\alpha_t\frac{\mu}{n}\right)\cr
 &\qquad\qquad\an{-(1-\eta_t)(1-\sigma^2)\left(1-\eta_t+\alpha_t\frac{\mu}{n}\right)}\cr
	&\le(1+\eta_t\dW)\left(1-\eta_t-2\left(L+\frac{2nL^2}{\mu}\right)\alpha_t\right),
    \end{align*}
    which is equivalent to
\begin{align*}
	&(1+\eta_t\dW -\an{(1-\eta_t)}(1-\sigma^2))\left(1-\eta_t+\alpha_t\frac{\mu}{n}\right)\cr
 &\le(1+\eta_t\dW)\left(1-\eta_t-2\left(L+\frac{2nL^2}{\mu}\right)\alpha_t\right).
 \end{align*}

        \begin{align*}
	&\frac{a_1}{b_1}=\frac{1-\eta_t+\mu\alpha_t/n}{1+\eta_t\|W-I\|^2}\cr
    &\qquad\le\frac{1-\eta_t-2\left(L+\frac{2nL^2}{\mu}\right)\alpha_t}{(1-\eta_t)\sigma^2 + \eta_t(1+\|W-I\|^2)}=\frac{a_2}{b_2}.
\end{align*}}
\medskip

\emph{2. More details on the inequality  $a_2 = \frac{1-\eta_t}{2} - \left(L+\frac{2nL^2}{\mu}\right)\alpha_t\ge1/8$ under the conditions $\alpha_t\le\frac{\mu}{4(L\mu+2nL^2)}$ and $\alpha_t\le\frac{\sqrt 3}{4L}$.}

As $\eta_t = \frac{1-\sqrt{1-4L^2\alpha_t^2}}{2}$,
$\frac{1-\eta_t}{2} - \left(L+\frac{2nL^2}{\mu}\right)\alpha_t\ge1/8$
if and only if 
\begin{align*}
	&2+2\sqrt{1-4L^2\alpha_t^2} - 8\left(L+\frac{2nL^2}{\mu}\right)\alpha_t\ge1\cr
	&\qquad\qquad\Updownarrow\cr
	&2 - 8\left(L+\frac{2nL^2}{\mu}\right)\alpha_t\ge1-2\sqrt{1-4L^2\alpha_t^2},
\end{align*}
which holds, if $2 - 8\left(L+\frac{2nL^2}{\mu}\right)\alpha_t\ge0$ and  $1-2\sqrt{1-4L^2\alpha_t^2}\le0$. 
The first inequality is guaranteed by $\alpha_t\le\frac{\mu}{4(L\mu+2nL^2)}$, whereas the second one is implied by $\alpha_t\le\frac{\sqrt 3}{4L}$.

\end{document}